\documentclass[10pt]{amsart}

\usepackage{amsmath, amsthm}
\usepackage{amssymb}
\usepackage{amscd}
\usepackage{latexsym}
\usepackage[latin1]{inputenc}

\usepackage[shortlabels]{enumitem}
\usepackage{extpfeil}

\usepackage{xypic}
\usepackage{graphicx}
\usepackage{changebar,color}

\usepackage{tikz-cd}

\usepackage{mathrsfs}

\newtheorem{theorem}{Theorem}[section]

\newtheorem{corollary}[theorem]{Corollary}

\newtheorem{proposition}[theorem]{Proposition}

\theoremstyle{remark}

\theoremstyle{definition}

\theoremstyle{definition}

\numberwithin{equation}{section}

\usetikzlibrary{decorations.markings}
\makeatletter
\tikzcdset{
open/.code={\tikzcdset{hook, circled};},
closed/.code={\tikzcdset{hook, slashed};},
open'/.code={\tikzcdset{hook', circled};},
closed'/.code={\tikzcdset{hook', slashed};},
circled/.code={\tikzcdset{markwith={\draw (0,0) circle (.375ex);}};},
slashed/.code={\tikzcdset{markwith={\draw[-] (-.4ex,-.4ex) -- (.4ex,.4ex);}};},
markwith/.code={
\pgfutil@ifundefined{tikz@library@decorations.markings@loaded}%
{\pgfutil@packageerror{tikz-cd}{You need to say %
\string\usetikzlibrary{decorations.markings} to use arrow with markings}{}}{}%
\pgfkeysalso{/tikz/postaction={/tikz/decorate,
/tikz/decoration={
markings,
mark = at position 0.5 with
{#1}}}}},
}
\makeatother

\renewcommand{\phi}{\varphi}

\def\N{\mathbb{N}}
\def\Z{\mathbb{Z}}
\def\PP{\mathbb{P}}
\def\Q{\mathbb{Q}}
\def\R{\mathbb{R}}

\def\C{\mathbb{C}}

\def\H{\mathbb{H}}
\def\V{\mathbb{V}}

\newcommand{\Gm}{{{\mathbb G}_m}}

\newcommand{\cF}{{\mathcal  F}}

\newcommand{\cGK}{{\mathcal  GK}}

\newcommand{\ccH}{{\mathcal  H}}
\newcommand{\cI}{{\mathcal  I}}

\newcommand{\cO}{{\mathcal O}}

\newcommand{\cV}{{\mathcal  V}}

\newcommand{{\OL}}{{{\mathcal O}_L}}

\newcommand{\Proj}{{\rm Proj\,}}
\newcommand{\Spec}{{\rm Spec\, }}

\newcommand{\rk}{{\rm rk }}

\newcommand{\ra}{\rightarrow}
\newcommand{\lrasim}{\stackrel{\sim}{\longrightarrow}}
\newcommand{\lra}{\longrightarrow}

\newcommand{\hlra}{{\lhook\joinrel\longrightarrow}}

\newcommand{\hgt}{{\mathrm{ht}}}

\newcommand{\Cc}{{\mathaccent23 C}}

\newcommand{\Omegac}{{\mathaccent23 \Omega}}

\newcommand{\Id}{{\mathrm{Id}}}

\renewcommand{\epsilon}{\varepsilon}

\newcommand{\GK}{\mathcal{GK}}

\title[Pencils of  hypersurfaces, Griffiths heights and geometric invariant theory. I]{Pencils of projective hypersurfaces, \\ Griffiths heights and geometric invariant theory. I}

\author{Thomas Mordant}

\begin{document}

\maketitle

\begin{abstract} We study the Griffiths heights associated to the middle-dimensional cohomology of pencils of projective hypersurfaces, by comparing them to heights defined by means of geometric invariant theory (GIT).  Kato and Koshikawa have conjectured a Northcott property for the Kato heights attached to motives over number fields, 
and investigated its consequences.  Bounding these Griffiths heights in terms of the GIT heights would constitute 
a geometric counterpart, valid over function fields of characteristic zero, of Kato and Koshikawa's conjecture.

Part  of our results follows from our earlier works \cite{Mordant22} on the computation of these Griffiths heights in the case of pencils with generic singularities, and \cite{Mordant23} on semistability criteria for singular projective hypersurfaces, combined with a general formalism of GIT heights over function fields. 

We also establish estimates between the Griffiths and GIT heights associated to pencils of projective hypersurfaces, which are valid beyond the case of generic singularities. To achieve this, we establish diverse results of independent interest. Notably we extend the computations in \cite{Mordant22} to pencils of projective hypersurfaces with semihomogeneous singularities, and we show the lower semicontinuity of the stable Griffiths height attached to polarized variations of Hodge structures over complex algebraic curves.

\end{abstract}

%
%
%

\tableofcontents

\section{Introduction}

\subsection{}\label{introduction Kato def} \addtocontents{toc}{\protect\setcounter{tocdepth}{1}}
In \cite{Kato14}, Kato has introduced a notion of height associated to pure motives over  number fields, that he investigated more thoroughly later (see \cite{Kato13} for an extension to mixed motives, and \cite{Kato18} and \cite{Kato20} for a detailed presentation; see also Koshikawa's work \cite{Koshikawa15} for variants of   and further results on Kato's height of pure motives). Kato's main inspiration was Faltings' height of abelian varieties over number fields, which played a key role in Faltings' proof of the Tate conjecture for abelian varieties and consequently of Mordell's conjecture (\cite{Faltings83}). 

Kato's height of pure motives, specialized to motives of weight one, essentially coincides with Faltings' height. Moreover Kato's definition for general weights is designed to satisfy a property of invariance under isogenies which generalizes the invariance under isogenies of Faltings' height. This invariance plays a central role in \cite{Faltings83}, and  is actually a number field analogue of a similar property of a suitable height associated by Zarhin to abelian varieties over function fields (\cite{Zarhin75}). 

However, while in the geometric case investigated by Zarhin, this invariance property is a direct consequence of the definition of the height, in the number field case considered by Faltings, its proof relies on the validity of the Weil conjecture for abelian varieties over finite fields and global class field theory.  Kato's proof of the invariance under isogenies of the height he attaches to motives over number fields relies on the general validity of the Weil conjecture as established by Deligne, and on the developments of $p$-adic  cohomology and comparison theorems.   

Kato and Koshikawa conjecture that the heights of pure motives they consider satisfy 
 finiteness properties \emph{\`a la} Northcott. According to a now-classical line of reasoning --- which goes back to Tate's seminal work \cite{Tate66} where Tate's conjecture is proved for abelian varieties over finite fields --- this conjectural finiteness of Kato's heights, combined with their invariance under isogenies, would admit striking consequences concerning  full faithfulness conjectures for realization functors from categories of ``geometric motives'' to systems of Galois representations. 
 
 Unfortunately, at this stage, nothing is known concerning the validity of these Northcott properties, except for motives of weight one, in which case it constitutes one of the most delicate technical points in Faltings' paper  \cite{Faltings83}. Detailed proofs of the Northcott property for Faltings heights may indeed be found in Moret-Bailly's expos\'e in Szpiro's seminar on Mordell's conjecture \cite{Moret-Bailly85a}, and  in the book  of Faltings and Chai \cite{Faltings-Chai90} on degenerations of abelian varieties (see notably \cite[Section V.4]{Faltings-Chai90}), and constitute some of the main achievements in the publications \cite{Szpiro85} and \cite{Faltings-Chai90}.\footnote{Both of these proofs 
 ultimately rely on the fundamental results concerning (semi-)abelian schemes established in Moret-Bailly's monograph \cite{Moret-Bailly85b}. This is equally the case of the proof in \cite{Bost96}, which also uses geometric invariant theory and higher dimensional Arakelov geometry.}

 \subsection{} 
  
 The Kato height $\hgt_K(M)$ of a motive $M$ over a number field $F$ is defined as the Arakelov degree of a certain one-dimensional $F$-vector space, which is constructed as the tensor product of powers of the determinant lines of the subquotients of the Hodge filtration on the de Rham realization $M_{\mathrm{d R}}$, and is equipped with metrics over the completions $F_v$ of $F$ at its various places $v$. These metrics are defined in terms of the $p$-adic realization of $M$ for $v$ of residue characteristic $p$, and of the Hodge realization of $M$ for $v$ an Archimedean place. 
 
 As pointed out by Kato in \cite{Kato18}, the specific combination of powers of subquotients of the Hodge filtration used to define this vector space is precisely the same as the one used by Griffiths in \cite{Griffiths70} to define the \emph{canonical line bundle}, which we will denote by $\GK_S(\V)$, of a variation of Hodge structures $\V$ over a (smooth) complex manifold $S$. 
 
 In retrospect, the Kato height of a motive $M$ over a number field $F$ appears as the arithmetic counterpart of the \emph{Griffiths height} of a variation of Hodge structures $\V$ over a connected smooth projective complex curve $C$, defined as the degree of its Griffiths canonical line bundle: 
  $$\hgt_{GK}(\V) := \deg_C \GK_C(\V).$$
 The analogy between the Kato height $\hgt_K(M)$ and  the Griffiths height $\hgt_{GK}(\V)$ is an instance of the classical analogy between number fields and function fields. With the above notation, the function field $\C(C)$ is the counterpart of the number field $F$, and a variation of Hodge structures $\V$ over $C$ --- or more generally over some  Zariski open subset $C \smallsetminus \Delta$ of $C$, with $\Delta$ a finite subset of $C$ --- plays the role of a motive $M$ over $F$.
 
 However it should be emphasized that  the Kato height and  the Griffiths canonical line bundle have been defined independently\footnote{When Kato first introduced his construction of heights, he was actually unaware of its relation to Griffiths' earlier work, which was pointed out by Koshikawa;  see \cite{Kato14} and \cite[Remark 1.6.4]{Kato18}.}, with completely different motivations.  
 
 While Kato's definition was inspired by Faltings' proof of the invariance under
 isogenies of the Faltings height of abelian varieties,  
 Griffiths was studying the properties of the Chern forms on the subquotients of the Hodge filtration induced by a polarization on some variation of Hodge structures~$\V$. Griffiths proved (\cite[Proposition (7.15)]{Griffiths70}) that when $\V$ is polarized, the associated Chern form on the line bundle $\GK_S(\V)$ is \emph{non-negative}, and vanishes if and only if the Hodge filtration of $\V$ is flat relatively to the natural connection on the underlying vector bundle $\cV$ of $\V$.

  \subsection{} The  construction of the Griffiths height $\hgt_{GK}(\V)$ was actually generalized by Peters (\cite{Peters84}) in the ``bad reduction'' case where the variation of Hodge structures $\V$ is only defined on the complement~$C \smallsetminus \Delta$ of a finite subset $\Delta$ of the curve $C$. 
  
  Peters used the theory of Deligne extensions of a vector bundle with connection to extend the line bundle $\GK_{C \smallsetminus \Delta}(\V)$ on $C \smallsetminus \Delta$ into a line bundle on $C$. This construction involves a choice of sign, and therefore there are two natural such line bundles extending $\GK_{C \smallsetminus \Delta}(\V)$: the \emph{upper} and \emph{lower} extensions, which we shall denote by $\GK_{C,+}(\V)$ and $\GK_{C,-}(\V)$ respectively. We shall moreover call their degrees on $C$ the \emph{upper} and \emph{lower Griffiths heights} of $\V$, and denote them by $\hgt_{GK,+}(\V)$ and~$\hgt_{GK,-}(\V)$.
 
 This construction allowed Peters to prove (see \cite[Theorem (4.1)]{Peters84}) that if the variation of Hodge structures $\V$ is polarized, then its \emph{upper} Griffiths height $\hgt_{GK,+}(\V)$ is non-negative, and it vanishes if and only if the Hodge filtration on $\V$ is flat and the local monodromy of $\V$ at every point of $\Delta$ is unipotent.

 Together with the upper and lower Griffiths heights, we may also consider the \emph{stable Griffiths height} $\hgt_{GK, stab}(\V)$, defined by pulling back the variation of Hodge structures $\V$ to a ramified covering of the curve $C$ such that the local monodromy of the resulting variation of Hodge structures is everywhere unipotent (see \cite[1.2]{Mordant22} and \ref{stGh} below). The definition of $\hgt_{GK, stab}(\V)$ is inspired by the definition of the stable Faltings height of abelian varieties. 
 
 The height $\hgt_{GK, stab}(\V)$ also satisfies a non-negativity result, and is compatible with base change under finite morphisms of smooth projective curves.

 The three Griffiths heights $\hgt_{GK, \pm}(\V)$ and $\hgt_{GK, stab}(\V)$ only depend on the restriction of $\V$ to an arbitrarily small non-empty  Zariski open subset of $C$. In this sense, they only depend  on the generic fiber 
of $\V$ over the function field $\C(C)$ of the curve $C$. This confirms their analogy with Kato heights of motives over number fields. This also allows  us to introduce the Griffiths heights attached to a variation of Hodge structures $\V_\eta$ on some unspecified non-empty Zariski open subset of $C$ --- these heights will be denoted by $\hgt_{GK, \pm}(\V_\eta)$ and $\hgt_{GK, stab}(\V_\eta)$, after the notation $\eta$ for the generic point of~$C$.

 \subsection{} \label{motivation GK} Curiously enough, although the geometric situation concerning Griffiths heights of variations of Hodge structures is supposed to be much more approachable than the arithmetic situation concerning Kato heights of motives, for a long time no significant result was apparently known about the former, with the exceptions of   the weight one case --- in which case the Griffiths height is basically the degree of the Hodge bundle of a pencil of abelian varieties --- and of Peters' positivity result.
 
 With the long-term motivation of getting some better understanding of Kato heights of motives, in \cite{Mordant22} we investigated the Griffiths heights associated to some significant instances of variations of Hodge structures of geometric origin and of arbitrary weight --- namely the middle-dimensional cohomology of pencils of projective hypersurfaces. 
 
More precisely, we computed the stable Griffiths height $\hgt_{GK, stab}\big(\H^{N-1}(H_\eta/C_\eta)\big)$ (along with its upper and lower variants) of the relative cohomology $\H^{N-1}(H_\eta/C_\eta)$ attached to a pencil of generically smooth hypersurfaces $H$ in a projective bundle over a smooth projective complex curve~$C$, defined by some vector bundle $E$ of rank~$N+1 > 1$:
$$\pi : \PP(E) := \Proj_C S^\bullet E^{\vee} \lra C.$$

Observe that the fibers of $H$ over $C$ are projective hypersurfaces of  dimension~$N-1,$ and have interesting cohomology groups only in middle dimension $N-1$, as follows from the weak Lefschetz theorem and Poincar\'e duality. 
Our computation was achieved when the pencil $H$ over $C$ had ``generic singularities'' in the following sense: the total space, of dimension $N$, of the hypersurface~$H$  was assumed to be smooth, and the flat morphism $\pi_{| H} : H \ra C$ to admit non-degenerate critical points only.\footnote{This last condition is equivalent to the requirement that the singularities of the  fibers of  $\pi_{| H}$ are ordinary double points.}

In this paper, we want to pursue our investigation of these Griffiths heights associated to pencils of projective hypersurfaces, and to explore their possible boundedness properties, which would be geometric analogues of the Northcott finiteness properties expected for  Kato heights. 

For instance, one would expect that, in the situation above, for a fixed curve $C$, and a fixed absolute complex dimension $N$ and relative degree $d$ of the pencil of hypersurfaces $H$, the pencils~$H$ over~$C$ with a given value of $\hgt_{GK, stab}\big(\H^{N-1}(H_\eta/C_\eta)\big)$ constitute a bounded family, as soon as $N \geq 2,$ $d\geq 3$ and $(N, d) \neq (3,3)$.\footnote{Here and in the next paragraphs, we exclude the case $d=2,$ dealing with pencils of quadrics, which is of limited interest since all smooth projective quadrics in $\PP^N_\C$ are projectively equivalent.}

The perspective of getting some ``bounds'' on the hypersurfaces themselves, using information on their middle-dimensional cohomology only, is comforted --- at least at some informal level --- by the validity of (some form of) the Torelli theorem in this context when assuming $N\geq 2,$ $d\geq 3$ and~$(N, d) \neq (3,3)$.\footnote{See \cite[Section 6.3.2]{Voisin03} and \cite{Voisin22} for a discussion of the classical results on (infinitesimal and generic) Torelli theorems for projective hypersurfaces, and some improvements of those.}

On the other hand, as in the arithmetic setting, it seems sensible for establishing such boundedness results to require some additional assumptions on the ``bad reduction'' of the pencil $H$ over $C$, of the same kind as our ``generic singularities'' assumption in our previous computation. These extra assumptions would allow for some control of the local monodromy of the variation of Hodge structures~$\H^{N-1}(H_\eta/C_\eta)$ around the points of $C$ where it is not defined.\footnote{Dealing with the ``bad reduction'' in Kato's arithmetic setting appears as an especially delicate issue. See \cite[Conjecture 4.2]{Kato14}, where Kato formulates a Northcott property concerning only motives over a number field with everywhere \emph{semistable} reduction, and Koshikawa's comments in \cite[Section 10.1]{Koshikawa15}. In a related vein, this issue motivated Koshikawa  to introduce some variant of Kato's original definition of the height attached to a motive over a number field.}

\subsection{}\label{introduction GIT def}  The approach to the boundedness properties of Griffiths heights advocated for in this paper consists in relating them to heights defined by means of geometric invariant theory (GIT), which satisfy suitable boundedness properties by construction.

When dealing with projective hypersurfaces, the construction of GIT heights may be described as follows. 
The reader will find more complete expositions, valid in a more general setting, in Sections \ref{GIT hyp} and \ref{sec:RelGIT}.

Let $(N,d)$ be a pair of integers such that $N \geq 2$ and $d\geq 3$.  A hypersurface $H$ of degree $d$ in $\PP^N_\C$ is defined by the vanishing of some non-zero homogeneous polynomial 
$$P \in \C[X_0, \dots, X_N]_d = \Gamma (\PP^N_\C, \cO(d)),$$
and $H$ is smooth if and only the discriminant $\mathrm{Disc}_{N,d}(P)$ of this polynomial does not vanish (see for instance \cite{GKZ08} and \cite{Demazure12}). Accordingly the set of smooth hypersurfaces of degree $d$ in $\PP_\C^N$ may be identified with the affine Zariski open subset $$\mathcal{H}^{\mathrm{sm}}_{N,d} := \PP(\C[X_0, \dots, X_N]_d) \setminus \Delta_{N,d},$$
of the projective space of forms of degree $d$ in $N+1$ variables\footnote{aka $(N+1)$-ary $d$-ics.}: $$\PP(\C[X_0, \dots, X_N]_d) \simeq \PP^{\binom{N+d}{N}-1}(\C),$$
defined as the complement of the discriminant divisor:
$$ \Delta_{N,d} := (\mathrm{Disc}_{N,d} = 0).$$

Two hypersurfaces $H_1 = (P_1 = 0)$ and $H_2 = (P_2 =0)$ in $\mathcal{H}^{\mathrm{sm}}_{N,d}$ are isomorphic as projective hypersurfaces\footnote{When $N\geq 4$ and $d\geq 3$, this holds if and and only if the complex varieties $H_1$ and $H_2$ are isomorphic.} when the points $[P_1]$ and $[P_2]$ belong to the same orbit under the natural action of $GL_{N+1}(\C)$, or of $SL_{N+1}(\C)$,  on $\PP(\C[X_0, \dots, X_N]_d).$ Accordingly the set of (projective) isomorphism classes of smooth hypersurfaces of degree $d$ in $\PP^N_\C$ may be identified with the quotient set~$\mathcal{H}^{\mathrm{sm}}_{N,d}/SL_{N+1}(\C).$

A classical theorem of Jordan \cite{Jordan80} asserts that the stabilizer in $SL_{N+1}(\C)$ of any element $[P]$ of $\mathcal{H}^{\mathrm{sm}}_{N,d}$ is finite. This implies that the orbits of $SL_{N+1}(\C)$ acting on  $\mathcal{H}^{\mathrm{sm}}_{N,d}$ are closed in the Zariski topology (see for instance \cite[Proposition 4.2]{MumfordFogartyKirwan94}). 

Geometric invariant theory, as developed by Mumford (\cite[Chapter 4, \S 2]{MumfordFogartyKirwan94}), may therefore be used to equip the set $\mathcal{H}^{\mathrm{sm}}_{N,d}/SL_{N+1}(\C)$ with a natural structure of complex affine variety, and actually to construct a natural projective compactification of this affine variety. 

Indeed, as indicated above, Jordan's theorem implies that $\mathcal{H}^{\mathrm{sm}}_{N,d}$ is a (clearly Zariski open and $SL_{N+1}(\C)$-invariant) subset of the set  $\PP(\C[X_0, \dots, X_N]_d)_{s}$ of properly stable points of the projective space $\PP(\C[X_0, \dots, X_N]_d)$ relatively to the action of $SL_{N+1, \C}$. 

The GIT quotient of $\PP(\C[X_0, \dots, X_N]_d)$ under this action is the projective variety defined by the graded algebra of $SL_{N+1, \C}$-invariant polynomials on $\C[X_0, \dots, X_N]_d$:
\begin{multline*}
\PP_\C^{\binom{N+d}{N}-1}/\!/ SL_{N+1, \C}  = \PP(\C[X_0, \dots, X_N]_d)/\!/ SL_{N+1, \C} \\ := \Proj \bigoplus_{n \in \N} \Gamma\big(\PP(\C[X_0, \dots, X_N]_d), \cO(n)\big)^{SL_{N+1, \C}}.
\end{multline*}
The subscheme $\PP(\C[X_0, \dots, X_N]_d)_{ss}$ of semistable  points of $\PP(\C[X_0, \dots, X_N]_d)$ is its $SL_{N+1, \C}$-invariant open subscheme defined as the complement of the base locus of these invariant polynomials. By construction, it is endowed with a canonical $SL_{N+1, \C}$-invariant morphism:
\begin{equation} \label{defqGIT}
q: \PP(\C[X_0, \dots, X_N]_d)_{ss} \lra \PP_\C^{\binom{N+d}{N}-1}/\!/ SL_{N+1, \C},
\end{equation}
and the $SL_{N+1, \C}$-equivariant line bundle $\cO(1)$ over $\PP(\C[X_0, \dots, X_N]_d)_{ss}$ descends into an ample $\Q$-line bundle $L$ over $\PP_\C^{\binom{N+d}{N}-1}/\!/ SL_{N+1, \C}$, 
 so that the following isomorphism holds:
$$q^\ast L \simeq \cO(1)_{\mid \PP(\C[X_0, \dots, X_N]_d)_{ss}}.$$

The subscheme  $\mathcal{H}^{\mathrm{sm}}_{N,d}$ is contained in $\PP(\C[X_0, \dots, X_N]_d)_{s}$, hence in $\PP(\C[X_0, \dots, X_N]_d)_{ss}$. Its image $q(\mathcal{H}^{\mathrm{sm}}_{N,d})$ is an affine open subscheme of the GIT quotient 
$\PP_\C^{\binom{N+d}{N}-1}/\!/ SL_{N+1, \C}$ 
and satisfies:
$$\mathcal{H}^{\mathrm{sm}}_{N,d} = q^{-1}\big( q(\mathcal{H}^{\mathrm{sm}}_{N,d})\big).$$
Moreover the fibers of the restriction:
$$q_{\mid \mathcal{H}^{\mathrm{sm}}_{N,d}}: \mathcal{H}^{\mathrm{sm}}_{N,d} \lra q(\mathcal{H}^{\mathrm{sm}}_{N,d})$$ are precisely the $SL_{N+1}(\C)$-orbits in $\mathcal{H}^{\mathrm{sm}}_{N,d}$.   

Accordingly the open subscheme $q(\mathcal{H}^{\mathrm{sm}}_{N,d})$ of the GIT quotient $\PP_\C^{\binom{N+d}{N}-1}/\!/ SL_{N+1, \C}$
defines the desired structure of affine algebraic variety on the  set~$\mathcal{H}^{\mathrm{sm}}_{N,d}/SL_{N+1}(\C)$ of isomorphism classes of smooth hypersurfaces of degree $d$ in $\PP^N_\C$.

The  compactification of $\mathcal{H}^{\mathrm{sm}}_{N,d}/SL_{N+1}(\C)$ provided by the projective variety $\PP_\C^{\binom{N+d}{N}-1}/\!/ SL_{N+1, \C}$
endowed with the ample $\Q$-line bundle $L$ allows us to attach an height $\hgt_{GIT}(H_\eta)$ to any smooth projective hypersurface $H_\eta$ of degree $d$ in the projective space $\PP^N_{\C(C)}$ over the function field $\C(C)$ of some connected smooth projective complex curve $C$. 

Namely such an hypersurface is defined by an equation 
$$P_\eta = 0,$$
for some homogeneous polynomial $P_\eta \in \C(C)[X_0, \dots, X_N]_d \smallsetminus \{0\}$, and $[P_\eta]$ defines a $\C(C)$-point of the open subscheme $\mathcal{H}^{\mathrm{sm}}_{N,d}$ of $ \PP(\C[X_0, \dots, X_N]_d)_{ss}$. Its image
$$Q_\eta := q([P_\eta])$$
by the morphism \eqref{defqGIT} defines a $\C(C)$-point of $\PP_\C^{\binom{N+d}{N}-1}/\!/ SL_{N+1, \C}$,
and by properness extends to a morphism of $\C$-varieties:
\begin{equation*}
Q: C \lra \PP_\C^{\binom{N+d}{N}-1}/\!/ SL_{N+1, \C}. 
\end{equation*}
We may therefore define the GIT height of $H_\eta$ as the degree of the curve $Q_\ast C$ with respect to the ample $\Q$-line bundle $L$:
$$\hgt_{GIT}(H_\eta) := \int_{\PP^{\binom{N+d}{N}-1}/\!/ SL_{N+1, \C}}[Q_\ast C] \cap c_1 (L) = \deg_C Q^\ast L.$$ 
 
By construction, it is a non-negative rational number, invariant under (projective) isomorphisms of the hypersurface $H_\eta$, and its construction is compatible with finite extensions of the function field~$\C(C)$.

 Recall that the stable Griffiths height $\hgt_{GK, stab}(\V_\eta)$ is compatible with finite extensions of the function field $\C(C)$, while in general $\hgt_{GK, +}(\V_\eta)$ and $\hgt_{GK, -}(\V_\eta)$ are not.
 In our investigation of the relations between Griffiths heights and GIT heights associated to pencils of projective hypersurfaces, it is therefore sensible to focus on the \emph{stable} Griffiths height $\hgt_{GK, stab}\big(\H^{N-1}(H_\eta/C_\eta)\big)$.

 \subsection{}\label{1.6} To state our first comparison result relating the stable Griffiths height $\hgt_{GK, stab}\big(\H^{N-1}(H_\eta/C_\eta)\big)$ and the GIT height $\hgt_{GIT}(H_\eta)$ associated to a smooth projective hypersurface $H_\eta$ over the function field $\C(C)$, 
we consider a projective bundle $\pi: \PP(E)\ra C$ defined by a vector bundle $E$ of rank~$N+1>1$ over the smooth projective complex curve $C$, and a hypersurface $H$ in $\PP(E)$ such that~$\pi_{\mid H}: H \ra C$ is a flat morphism, as in~\ref{motivation GK} above. 

To stay in the setting of \ref{introduction GIT def} above, we shall also assume that the dimension $N$ of $H$ satisfies~$N \geq 2$, and the  degree $d$ of  the fibers of $\pi_{\mid H}$ satisfies $d \geq 3$.

\begin{theorem}\label{GK GIT generic} With the above notation, if $H$ defines a pencil of hypersurfaces with generic singularities, namely if $H$ is smooth and if $\pi_{\mid H}$ admits non-degenerate critical points only, then the following equality holds:
\begin{equation} \label{eq GK GIT generic}
\hgt_{GK, stab}\big(\H^{N-1}(H_\eta/C_\eta)\big) = F_{stab}(d,N)\,  \hgt_{GIT}(H_\eta),
\end{equation}
where  $F_{stab}(d,N)$ denotes the rational number in $(1/12) \Z$ defined, when $N$ is odd, by:
\begin{equation}\label{defFstabodd}
F_{stab}(d,N) := \frac{N+1}{24 d^2} \left[ (d-1)^N  (d^2 N - d^2 - 2 d N - 2 )+ 2 (d^2-1) \right],
\end{equation}
and, when $N$ is even, by:
\begin{equation}\label{defFstabeven}
F_{stab}(d,N) := \frac{N+1}{24 d^2} \left [ (d-1)^N  (d^2 N + 2 d^2 - 2 d N - 2) - 2 (d^2-1) \right ]. \end{equation}
\end{theorem}

Under the above assumptions on $N$ and $d$, the constant  $F_{stab}(d,N)$ is easily seen to be positive unless $(N,d) =(3,3)$, in which case it vanishes.\footnote{Actually, when $(N,d) =(3,3)$, the variation of Hodge structures $\H^{N-1}(H_\eta/C_\eta)$ has rank 7, is concentrated in bidegree $(1,1)$, and may be described in terms of the lines in the cubic surface $H_\eta$. This explains the isotriviality of~$\H^{N-1}(H_\eta/C_\eta)$ and the failure of the infinitesimal Torelli theorem in this case.} The equality \eqref{eq GK GIT generic} shows that, under the assumptions of Theorem \ref{GK GIT generic}, for a fixed pair of integers $(N,d)\neq (3,3)$, the pencils $H$ over $C$ for which $\hgt_{GK, stab}\big(\H^{N-1}(H_\eta/C_\eta)\big)$ is bounded  satisfy a boundedness property. This geometric result supports the expectation concerning the Northcott property of the Kato heights mentioned in \ref{introduction Kato def} above.   

The simplicity of the proportionality relation \eqref{eq GK GIT generic} conceals the fact that it is established as the  conjunction of two results of different nature --- namely the computations, under the assumptions of Theorem \ref{GK GIT generic}, of both the stable Griffiths height $\hgt_{GK, stab}\big(\H^{N-1}(H_\eta/C_\eta)\big)$ and the GIT height~$\hgt_{GIT}(H_\eta)$ in terms of a third invariant attached to the pencil $H/C$: the invariant $\hgt_{int}(H/C)$ already considered in \cite{Mordant22}.

Let us briefly recall the definition of $\hgt_{int}(H/C)$. Observe that, at a technical level,  this definition is much simpler  than the ones of  $\hgt_{GK, stab}\big(\H^{N-1}(H_\eta/C_\eta)\big)$ and $\hgt_{GIT}(H_\eta)$. However, in general,~$\hgt_{int}(H/C)$ actually depends on $H$ as a scheme fibered over $C$, and not only on its generic fiber.

The invariant $\hgt_{int}(H/C)$ may be attached  to any effective divisor  $H$ in a projective bundle
$\pi : \PP(E) \ra C$
of relative dimension $N$ over a connected smooth projective complex curve $C$. It is defined in terms of the relatively ample $\Q$-line bundle on $\PP(E)$:
\begin{equation} \label{def O1 *}
\cO_{\PP(E)}(1)^\ast := \cO_E(1) \otimes \pi^\ast (\det E)^{\otimes 1/(N+1)} \simeq \omega_{\PP(E)/C}^{\otimes -1/(N+1)},
\end{equation}
 --- where $\omega_{\PP(E)/C}$ denotes the relative dualizing sheaf of $\PP(E)$ over $C$ --- as the (rational) intersection number: 
$$\hgt_{int}(H/C) := \int_{\PP(E)} c_1(\cO_{\PP(E)}(1)^\ast)^N \cap [H].$$

Observe that  the isomorphism in \eqref{def O1 *} shows that the $\Q$-line bundle $\cO_{\PP(E)}(1)^\ast$ only depends on the $C$-scheme $\PP(E)$ and not on the choice of the underlying vector bundle $E$. Consequently the rational number $\hgt_{int}(H/C)$ only depends on the $C$-schemes $\PP(E)$ and $H$.

The comparison \eqref{eq GK GIT generic} between the stable Griffiths height and the GIT height arises as a consequence of two independent results. 

The first one is the computation in \cite{Mordant22} --- using notably Steenbrink's theory and the Gro\-then\-dieck-Riemann-Roch theorem --- of the stable Griffiths height in terms of  
$\hgt_{int}(H/C)$ in the case of ``generic singularities'', namely the following equality:
\begin{equation}\label{eq GK int generic}
\hgt_{GK, stab}\big (\H^{N-1}(H_\eta/C_\eta)\big ) = F_{stab}(d,N)\,  \hgt_{int}(H/C);
\end{equation}
 see \cite[Theorem 1.4.2]{Mordant22}.

The second result is a comparison of $\hgt_{int}(H/C)$ with the GIT height, which constitutes one of the main results in this paper. We show that, under the assumptions of Theorem \ref{GK GIT generic}, the following equality holds: 
\begin{equation}\label{eq GIT int}
\hgt_{GIT}(H_\eta) = \hgt_{int}(H/C).
\end{equation}

Actually we show in Sections \ref{GIT hyp} and \ref{sec:RelGIT}  that the following
 inequality holds for any pencil of hypersurfaces $H$ with semistable (e.g. non-singular) generic fiber $H_\eta$:
\begin{equation}\label{ineq GIT int}
\hgt_{GIT}(H_\eta) \leq \hgt_{int}(H/C),
\end{equation}
and that equality holds in \eqref{ineq GIT int} if and only if every fiber of the pencil $H$ over $C$ is a semistable hypersurface (see Theorem \ref{theorem hGITH} below). In the situation of Theorem \ref{GK GIT generic}, this semistability assumption is a consequence of criteria for the GIT
 semistability of singular projective hypersurfaces of arbitrary  dimension $N-1$  and degree $d$, established  in \cite{Lee08} when $d \geq N+1$  and in \cite{Mordant23} for any $d \geq 3$ (including the Fano case $d \leq N$).
 
 \subsection{}\label{1.7} In this paper, we also investigate the relation between $\hgt_{GK, stab}\big(\H^{N-1}(H_\eta/C_\eta)\big)$ and $\hgt_{GIT}(H_\eta)$ beyond the situation of pencils with ``generic singularities'' covered in Theorem  \ref{GK GIT generic}. 
 
 To achieve this, we extend the computation of $\hgt_{GK, stab}\big(\H^{N-1}(H_\eta/C_\eta)\big)$ in \cite{Mordant22} to pencils of hypersurfaces whose singular fibers are significantly more general than hypersurfaces with ordinary double points, namely hypersurfaces with semihomogeneous singularities of arbitrary multiplicity.  We are able to establish a generalization of the equality \eqref{eq GK int generic}, where the right-hand side contains, in addition to a multiple of $\hgt_{int}(H/C)$, a sum of local terms taking the multiplicities of the singularities of the fibers of $\pi_{\mid H}: H \ra C$ into account. 
 
 We refer the reader to 
 Theorem \ref{intro GK hyp P(E) hom crit} for the precise statement of this generalization.  Its derivation follows the strategy developed in \cite{Mordant22} to achieve the computation of $\hgt_{GK, stab}\big(\H^{N-1}(H_\eta/C_\eta)\big)$ for pencils with generic singularities, but also introduces some significant variants. These allow us to circumvent a major issue in this general strategy: the determination of the elementary exponents of the variation of Hodge structures $\H^{N-1}(H_\eta/C_\eta)$ at the singular fibers of $\pi_{\mid H}$. 
 
 Indeed, when these fibers have semihomogeneous singularities, it is possible to introduce a finite covering $C'$ of $C$, of generic point $\eta'$, on which the pull-back $\H^{N-1}(H_{\eta'}/C'_{\eta'})$ of $\H^{N-1}(H_\eta/C_\eta)$  admits unipotent local monodromy, and to describe the Deligne extension of $\H^{N-1}(H_{\eta'}/C'_{\eta'})$ by constructing a resolution~$\widetilde{H}_{C'}$ of $H_{C'} := H \times_C C'$ such that Steenbrink's theory may be applied to the pencil $\widetilde{H}_{C'}$ over~$C'$. 
 
 The details of this construction are somewhat technical, and we defer the proof of Theorem  \ref{intro GK hyp P(E) hom crit} to the second part of this paper \cite{MordantGIT2}. 
 
 Our computation of $\hgt_{GK, stab}\big(\H^{N-1}(H_\eta/C_\eta)\big)$ for pencils of hypersurfaces with semihomogeneous singularities, together with our results on GIT heights   in Section \ref{GIT hyp}, allow us to construct examples of smooth hypersurfaces $H_\eta$ over the function field $\C(C)$ for which the proportionality relation    \eqref{eq GK GIT generic} no longer holds. 
 
 However, in spite of the complexity of the final expression for  $\hgt_{GK, stab}\big(\H^{N-1}(H_\eta/C_\eta)\big)$ in this setting, we are able to prove that the following estimate still holds, under suitable assumptions on the dimension $N$ and the relative degree $d$ of the pencil and the maximal multiplicity $\delta$ of the singularities of its fibers\footnote{Namely when $N \geq 2$ and $d \geq 3 (\delta-1),$ or when $d \geq N+1$ and $d \geq \delta \, (1 + 1/N)$. These numerical conditions are required to ensure the semistability of the singular fibers of the pencil, or equivalently  the equality \eqref{eq GIT int}.}:
\begin{equation}\label{ineq GK GIT semihom}
\hgt_{GK, stab}\big(\H^{N-1}(H_\eta/C_\eta)\big) \leq F_{stab}(d,N)\, \hgt_{GIT}(H_\eta)
\end{equation}
(see Corollary \ref{pref intro GK hyp P(E) hom crit GIT} below).

\subsection{}\label{1.8} These results, based on the computation of $\hgt_{GK, stab}\big(\H^{N-1}(H_\eta/C_\eta)\big)$ for pencils of hypersurfaces with semihomogeneous singularities, naturally lead to ask oneself (i) at which level of generality the estimate \eqref{ineq GK GIT semihom} may be expected to hold, and (ii) whether some estimate in the opposite direction, of the form 
\begin{equation}\label{ineq GK GIT ??}
\hgt_{GK, stab}\big(\H^{N-1}(H_\eta/C_\eta)\big) \geq F_{??}(d,N)\, \hgt_{GIT}(H_\eta),
\end{equation}
for some positive constant $F_{??}(d,N)$, would hold for a smooth projective hypersurface $H_\eta$ of dimension $N-1$ and degree $d$ over $\C(C),$ under some suitable geometric condition on its ``fibers of bad reduction''. 

A positive answer to (ii) would provide some additional support for the validity of the Northcott property for Griffiths and Kato heights. Unfortunately we have not been able to reach a neat conclusion on this question yet.

Concerning (i), at this stage, it seems sensible to expect that indeed the estimate \eqref{ineq GK GIT semihom} holds for an arbitrary smooth projective hypersurface $H_\eta$ of dimension $N-1$ and degree $d$ over $\C(C).$ Indeed, besides its validity for suitable pencils with semihomogeneous singularities, in Sections \ref{sec: semicontinuity GK}
and  \ref{sec: ineq GK semi-hom}, we shall establish it when $H_\eta$ admits as a ``model'' over $C$ a divisor $H$ in some projective bundle $\PP(E)$, when the fibers of $H$ over $C$ are semistable projective hypersurfaces and  the height $\hgt_{int}(H/C)$ is large enough  compared to the ``instability'' $\mu_{\max}(E) - \mu(E)$ of the vector bundle $E$. 

The proof of \eqref{ineq GK GIT semihom}  in this case will combine the equality \eqref{eq GK int generic} valid for pencils with ``generic singularities'', and a degeneration argument which will notably use a general lower semicontinuity property, of independent interest, of the stable Griffiths height of variations of Hodge structures on a complex projective curve.

\subsection{} Finally we would like to emphasize the exploratory character of this work. 

When writing this paper, our primary goal has been to gather some  material relevant to the understanding of the   Griffiths and Kato heights and their Northcott properties, concerning significant classes of   variations of Hodge structures or motives of weight $>1.$ According to the general principle claiming that geometry over function fields is simpler than over number fields, and therefore has to be studied first, we  focused on the former --- that is, on Griffiths heights. 

It is rather remarkable that this limited goal leads to various basic questions in complex algebraic and analytic geometry that apparently had not been investigated in the literature --- for example the existence of criteria for the semistability of singular projective hypersurfaces of arbitrary dimension and degree, as the ones established in \cite{Mordant23}; or the lower semicontinuity property of the stable Griffiths height established in Section \ref{sec: semicontinuity GK}. 

Actually, to keep this paper within a limited number of pages, several of these basic questions have been incompletely investigated. This is for instance the case for the construction of heights over function fields using geometric invariant theory and their relations to boundedness properties. A more satisfactory approach to these problems  would involve the formalism of good moduli spaces associated to algebraic stacks, as developed by  Alper, Halpern-Leistner and Heinloth \cite{Alper13, AHLH23}. We plan to return to this topic in some future work.

\subsection{Acknowledgements} I am very grateful to  Isabelle Mordant for her remarks on the introduction, and to  Damien Simon for his \LaTeX-nical assistance during the preparation of this paper. 

This paper explores the properties of the Griffiths heights associated to pencils of projective hypersurfaces, and notably in its second part \cite{MordantGIT2}, may be seen as  a sequel of my PhD memoir \cite{Mordant22}, written under the supervision of Jean-Beno\^it Bost. I thank him for discussions on the themes explored in this paper, and his advice during the preparation of this paper.

I also thank Gerard Freixas for his helpful comments on the results in this paper. Notably the possibility of using my earlier result \eqref{eq GK int generic} established  in \cite{Mordant22}, valid for pencils of hypersurfaces with generic singularities, together with a suitable semicontinuity property of the stable Griffiths height,  to derive an inequality of the form \eqref{ineq GK GIT semihom} under much more general assumptions, arose in a conversation with Gerard. I am very grateful to him for this suggestion.

During the work leading to this paper, I benefitted from the support of the Fondation Math\'e\-ma\-ti\-que Jacques Hadamard, as a \emph{Hadamard Lecturer}.

\subsection{} This paper is organized as follows.

 In its first half (Sections 2 to 5), we present in full detail the results alluded to in Subsections \ref{1.6} to \ref{1.8} above, together with the main steps of their proofs. Its second half (Sections 6 to 8) is devoted to the proofs of various technical results used in their proofs. Some of these are of independent interest, from the perspective of geometric invariant theory (Section 6), or of the study of variations of Hodge structures (Section 7). Section 8 is of a more technical nature. 

The derivation of our main result about the stable Griffiths heights associated to pencils of hypersurfaces with semihomogeneous singularities (Theorem \ref{intro GK hyp P(E) hom crit}) is deferred to the second part of this paper \cite{MordantGIT2}. This derivation heavily relies on the strategy developed in \cite{Mordant22}, while the present paper does not require any familiarity with it.   

In Sections \ref{sec: GriffitsHeight} and \ref{GIT hyp}, we  introduce the ``main characters'': the Griffiths height and the GIT height attached to pencils of projective hypersurfaces with non-singular generic fibers. 

Section \ref{sec: GriffitsHeight} recalls  the main definitions and basic facts concerning Griffiths heights of variations of Hodge structures on curves, with a special emphasis on stable Griffiths heights. 

In Section \ref{GIT hyp}, we introduce the GIT height in a more general framework than the one of projective hypersurfaces. Indeed we consider not only  the action of the linear group $GL_{N+1, \C}$ on the space $\C[X_0, \dots, X_N]_d$ of homogeneous forms of degree $d$ by change of coordinates --- which underlies the classification of projective hypersurfaces of degree $d$ in $\PP^N_\C$ --- but  an arbitrary homogeneous representation of the linear group.\footnote{In the arithmetic setting, similar constructions of heights using geometric invariant theory have been studied in the papers \cite{Bost94,  Zhang96, Bost96, Gasbarri00, Maculan17}.} 
We also relate this GIT height to 
  ``elementary'' heights attached to sections of projective bundles over a curve. In the case of pencils of projective hypersurfaces,   these elementary heights reduce to 
  the invariant $\hgt_{int}(H/C)$ introduced above.

In Section \ref{sec: GKGIT for double pt}, we establish Theorem \ref{GK GIT generic}. Our proof relies on the main results in \cite{Mordant22} and \cite{Mordant23}, combined with our results on GIT heights presented in Section \ref{GIT hyp}.

Section \ref{sec: GKGIT for semi hom} is devoted to the relations between
$\hgt_{GK, stab}\big(\H^{N-1}(H_\eta/C_\eta)\big)$ and $\hgt_{GIT}(H_\eta)$ beyond the situation of pencils with ``generic singularities'' discussed in \ref{1.7} and \ref{1.8} above.  
Notably we consider  a pencil of projective hypersurfaces with semihomogeneous singularities $H/C$ in some projective bundle $\PP(E)$ over $C$, and in Theorem \ref{intro GK hyp P(E) hom crit} we give a formula for the  Griffiths height $\hgt_{GK, stab}\big(\H^{N-1}(H_\eta/C_\eta)\big)$ in terms of $\hgt_{int}(H/C)$ and of the multiplicities of the singular points of its fibers.  This formula, whose proof is the content of \cite{MordantGIT2},  will allow us to establish that the estimate \eqref{ineq GK GIT semihom} holds for these pencils. 

We also explain how the estimate \eqref{ineq GK GIT semihom} still holds under a simple numerical condition on the instability $\mu_{\max}(E) - \mu(E)$ of the vector bundle $E$. Its proof under this condition will use a specialization argument, to  derive an estimate from the equality \eqref{eq GK int generic} which  holds under some genericity assumption only. 

As discussed above, the last three sections are of a different nature.

Section \ref{sec:RelGIT} contains  a detailed construction of the GIT height associated to some homogeneous representation of the linear group, and the proof of its basic properties. Notably we compare the GIT height and the 
``elementary'' heights attached to sections of projective bundles over a curve, alluded to above. Applied to the situation of hypersurfaces, this establishes the inequality \eqref{ineq GIT int} and the subsequent analysis of the equality case. 

Section \ref{sec: semicontinuity GK} is devoted to the proof of the lower semicontinuity property of the stable Griffiths height, needed in the proof of the estimate \eqref{ineq GK GIT semihom} by a specialization argument. 
Observe that our proof of this lower semicontinuity property, which is basically an algebraic fact, is established by using some analytic tools. Ultimately it appears as a consequence of the non-negativity of the  Chern form of the Griffiths line bundle attached to a polarized variation of Hodge structures, of the expression of the stable Griffiths height as an integral of this Chern form  --- these are direct consequences of the classical work of Griffiths and Peters \cite{Griffiths70, Peters84} --- and of (an avatar of) Fatou's lemma in integration theory.   

Sections \ref{sec:RelGIT} and \ref{sec: semicontinuity GK}  may be read independently of the other parts of the paper. 

Finally,  in Section  \ref{sec: ineq GK semi-hom}, we tie some loose ends in the proofs of the various cases of estimate \eqref{ineq GK GIT semihom} stated in Section \ref{sec: GKGIT for semi hom}. We provide the details of the computations needed to derive this estimate for pencils of hypersurfaces with semihomogeneous singularities from the cumbersome formula in Theorem \ref{intro GK hyp P(E) hom crit} for the attached stable Griffiths height. We conclude by spelling out the specialization argument establishing the validity of \eqref{ineq GK GIT semihom} under some  numerical condition on the instability $\mu_{\max}(E) - \mu(E)$. Besides the results in Sections  \ref{sec: GKGIT for double pt} and \ref{sec: semicontinuity GK}, this argument also uses some (not so) well-known constructions of pencils of hypersurfaces with generic singularities presented in Appendix A of~\cite{MordantMem23}.

\addtocontents{toc}{\protect\setcounter{tocdepth}{2}}

\section{The stable Griffiths height of a polarized variation of Hodge structures over a complex algebraic curve}\label{sec: GriffitsHeight}

\subsection{The Griffiths line bundle of a variation of Hodge structures and Peters' construction} 

\subsubsection{} Let $S$ be a smooth complex manifold, $\cV$ a vector bundle on $S$ and 
$$F^\bullet := F^0:= \mathcal{V} \supseteq F^1 \supseteq \dots \supseteq F^n \supseteq F^{n+1} = 0$$
a decreasing filtration by subbundles. As in \cite{Griffiths70}, we can attach to $(\cV, F^\bullet)$ its \emph{Griffiths line bundle}:
$$\mathcal{\cGK}_S(\cV, F^\bullet) := \bigotimes_{i=1}^n \det F^i \simeq \bigotimes_{r=0}^n (\det F^r/F^{r+1})^{\otimes r}.$$

If $C$ is a connected smooth projective complex curve, the \emph{Griffiths height} $\hgt_{GK}(\cV, F^\bullet)$ is simply the degree of this line bundle.

Using Deligne's work on extensions of analytic vector bundles with connections admitting only regular singular points, Peters showed how to extend this definition to integral variations of Hodge structures defined on non-empty Zariski open subsets of a connected smooth projective complex curve.

\subsubsection{Peters' construction: the unipotent case} Let $C$ be a connected smooth projective complex curve with generic point $\eta$, $\Delta$ a finite subset of $C,$ and $\Cc := C \smallsetminus \Delta$ its complement. Let also $\V = (V_\Z, \cF^\bullet)$ be an integral variation of Hodge structures (VHS) defined on $\Cc,$ where $V_\Z$ denotes a local system of finitely generated free $\Z$-modules over $\Cc$, and  $\cF^\bullet$  the Hodge filtration on the associated analytic vector bundle:
$$\cV := V_\Z \otimes_{\Z_{\Cc}} \cO^{\mathrm{an}}_{\Cc}.$$

Let $\nabla$ be the connection on the vector bundle $\cV$ defined by the local system $V_\Z.$
First, in the case where the local monodromy $T_x$ of the connection $\nabla$ at every point $x$ of $\Delta$ is a unipotent endomorphism, the vector bundle with connection $(\cV, \nabla)$ on $\Cc$ admits a canonical Deligne extension~$(\overline{\cV}, \overline{\nabla})$ to $C$, where:
$$ \overline{\nabla} : \overline{\cV} \lra  \overline{\cV} \otimes_{\cO^{\mathrm{an}}_{C}} \Omega^1_{C} (\Delta)$$
is a logarithmic connection. This extension is characterized by the property that, for every point $x$ in $\Delta,$ the residue endomorphism of the logarithmic connection $\overline{\nabla}$ at $x$ is nilpotent.

According to a result of Griffiths \cite[Theorem (4.13), a)]{Schmid73}, the subbundles $(\cF^p)_{p \geq 0}$ of the vector bundle $\cV$ on $\Cc$ are algebraic with respect to the algebraic structure defined by the vector bundle $\overline{\cV}$ on $C.$ Consequently they admit extensions $(\overline{\cF}^p)_{p \geq 0}$ as subbundles of $\overline{\cV}$.

The $\emph{Griffiths line bundle}$ of the variation of Hodge structures $\V$ is the line bundle on $C$ defined by:
$$\GK_C(\V_\eta) := \GK_C(\overline{\cV}, \overline{\cF}^\bullet),$$
and the \emph{Griffiths height} $\hgt_{GK}(\V_\eta)$ of $\V$ is the degree of this line bundle:
$$\hgt_{GK}(\V_\eta) := \deg_C \GK_C(\V_\eta).$$

This line bundle and this height only depend on the restriction of the variation of Hodge structures 
to any non-empty Zariski open subset of $C$, or equivalently to any open neighborhood in $C$ of the generic point $\eta$. This justifies the notation. 

The construction of  $\GK_C(\V_\eta)$ and of $\hgt_{GK}(\V_\eta)$ is also compatible with base change in the following sense. If $C'$ is a connected smooth projective complex curve with generic point $\eta',$ and if
$\sigma: C' \ra C$
is a finite morphism, then there is a canonical isomorphism of line bundles on $C'$:
$$\GK_{C'}\big((\sigma^* \V)_{\eta'}\big) \lrasim \sigma^*\GK_C(\V_\eta),$$
and consequently the following equality of integers holds:
\begin{equation}\label{intro GK sigma V}
 \hgt_{GK}\big((\sigma^* \V)_{\eta'} \big) = \deg(\sigma) \, \hgt_{GK}(\V_\eta).
 \end{equation}
 
 \subsubsection{Peters' construction: the general case}
 
When the local monodromy is not unipotent, the construction of Deligne extensions of vector bundles with connections depends on the determination of a logarithm (from $\C^\ast$ to $\C$).

There are two natural such choices,\footnote{They are defined by the determinations of the complex logarithm with values in $\R + i [0, 2\pi)$ and in $\R + i (-2\pi, 0]$; see for instance \cite[1.2]{Mordant22}.}  leading to two definitions of Deligne extensions, which coincide with the canonical Deligne extension when the local monodromy is unipotent: the \emph{upper extension} $\overline{\cV}_{+}$ and the \emph{lower extension} $\overline{\cV}_-$.  These complex analytic vector bundles over $C$ extend the vector bundle $\cV$ on $\Cc$ and have the same meromorphic structure at $\Delta$, and consequently define a canonical algebraic structure on $\cV$ over the complex algebraic curve $\Cc$. 

According to the result of Griffiths cited above, the subbundles $(\cF^p)_{p \geq 0}$ of the vector bundle~$\cV$ on $\Cc$ are again algebraic with respect to this algebraic structure  and consequently  admit extensions~$(\overline{\cF}^p_{\pm})_{p \geq 0}$ as subbundles of $\overline{\cV}_{\pm}$.

This leads to the definition of the upper and lower Griffiths line bundles:
$$\GK_{C, \pm}(\V_\eta) := \GK_C(\overline{\cV}_{\pm}, \overline{\cF}^\bullet_{\pm}),$$
and the upper and lower Griffiths heights:
$$\hgt_{GK, \pm}(\V_\eta) := \deg_C \GK_{C, \pm}(\V_\eta).$$

Peters specifically studied the upper variant in \cite{Peters84}. Extending an earlier result of Griffiths in the ``good reduction'' case, he showed in \cite[Theorem (4.1)]{Peters84} that the upper Griffiths height~$\hgt_{GK, +}(\V_\eta)$ of a \emph{polarized} VHS $\V$ is always non-negative, and vanishes if and only if all the subbundles $(\cF^p)_{p \geq 0}$ are flat relatively to the connection $\nabla$ and the monodromy of $\nabla$ at every point of $\Delta$ is unipotent.

In contrast, it follows from Steenbrink's theory (see \cite{Steenbrink76, Steenbrink77}) that in the geometric situation where $\V$ is the relative cohomology of a smooth pencil $Y$ of projective varieties over the algebraic curve $\Cc$, the \emph{lower} Deligne extension $\overline{\cV}_-$ and the attached extensions $(\overline{\cF}^p_{-})_{p \geq 0}$ of the subbundles~$(\cF^p)_{p \geq 0}$ are canonically isomorphic  to the relative logarithmic cohomology of a compactification $\overline{Y}$ of $Y$ over $C$ whose singular fibers are divisors with strict normal crossings, equipped with its natural Hodge filtration defined in terms of the relative logarithmic de Rham complex.

\subsection{The stable Griffiths height}\label{stGh} In general, the construction of the Deligne extensions, and consequently of the upper and lower  Griffiths heights, 
 is not compatible  with base change under finite morphisms of curves. In other words, the equality \eqref{intro GK sigma V} does not hold in general with $ \hgt_{GK, +}$ or $ \hgt_{GK, -}$ instead of $\hgt_{GK}$ when the local monodromy at $\Delta$ is not unipotent.

To fix this issue, as in \cite{Mordant22}, one can introduce 
the so-called \emph{stable Griffiths height} $\hgt_{GK, stab}(\V_\eta),$ the definition of which is inspired by the definition of the stable Faltings height of abelian varieties.

To define $\hgt_{GK, stab}(\V_\eta),$ recall that, 
according to a result of Borel \cite[Lemma (4.5)]{Schmid73}, the local monodromy $T_x$ of the connection $\nabla$ at every point $x$ of $\Delta$ is a quasi-unipotent endomorphism, namely there exists an integer $r_x \geq 1$ such that the endomorphism $T_x^{r_x}$ is unipotent. Let $C'$ be a connected smooth projective complex curve and
$\sigma: C' \ra C$
 a finite morphism that is ramified at every point~$x'$ of $\Delta' := \sigma^{-1}(\Delta)$ with ramification index a multiple of $r_{\sigma(x')}.$ Such a ramified covering~$C'$ of $C$ may be constructed as a cyclic covering of $C$. Moreover any two such ramified covering $C'_1$ and $C'_2$ are dominated by a third one $C''$.  
 
According to the ramification assumption, the local monodromy at every point of $\Delta'$ of the pulled back vector bundle with connection $(\sigma^* \cV, \sigma^*\cF^\bullet)$ on $\sigma^{-1}(\Cc)$ is unipotent. Consequently it admits a canonical Deligne extension, and as before, we can define the Griffiths line bundle $\GK_{C'}\big((\sigma^* \V)_{\eta'}\big)$ and the Griffiths height $\hgt_{GK}\big((\sigma^* \V)_{\eta'}\big)$.

The \emph{stable Griffiths height} of the variation of Hodge structures $\V$ is the rational number defined by:
$$\hgt_{GK, stab}(\V_\eta) := \deg(\sigma)^{-1}\,   \hgt_{GK}\big((\sigma^* \V)_{\eta'}\big) = \deg(\sigma)^{-1}\, \deg_{C'}\GK_{C'}\big((\sigma^* \V)_{\eta'}\big).$$

This rational number indeed does not depend on the choice of the finite covering $\sigma$: this follows from the compatibility \eqref{intro GK sigma V} of the Griffiths height with finite coverings in the unipotent case.  Furthermore, like the definition of $\hgt_{GK}(\V_\eta)$ in the unipotent case, it only depends on the restriction of the variation of Hodge structures to any non-empty Zariski open subset of $C$,  which justifies the notation. 

By construction, the height $\hgt_{GK, stab}$ is  compatible with base change, in the sense that \eqref{intro GK sigma V} holds  for every finite  covering $C'$ of $C$, with $\hgt_{GK, stab}$ instead of $\hgt_{GK}$, as an equality of rational numbers.

\section{The GIT height of a pencil of hypersurfaces} \label{GIT hyp}

\subsection{The GIT height of a $k(C)$-point of the semistable locus $\PP^{v-1}_{k, ss}$ attached to a homogeneous representation $GL_{e,k} \ra  GL_{v,k}$} 
\label{GIT height quot}

 In this subsection, we consider a connected smooth projective curve $C$ over an algebraically closed field $k$, two positive integers $e$ and $v$, and a morphism of algebraic groups over $k$:
 $$\rho : GL_{e,k} \lra  GL_{v,k}.$$
 We assume that the representation $\rho$ is homogeneous, namely maps $\Gm \Id_e$ into $\Gm \Id_v$, and we denote by $\gamma$ its weight, namely the integer such that:
 $$\rho (t\, \Id_e) = t^\gamma\, \Id_v$$
 for every $t \in k^\ast.$
 
 With this notation, we shall associate to every point of the semistable locus  $\PP^{v-1}_{k, ss}$ over the function field $k(C)$ of $C$ a canonical height defined in terms of the GIT quotient  $\PP^{v-1}_k /\!/ SL_{e,k}$. The proofs of the results of this subsection are deferred to Section \ref{sec:RelGIT}. 
 
\subsubsection{}  \label{relative ss quot}
Geometric invariant theory allows us to consider the Zariski open subscheme 
$\PP^{v-1}_{k, ss}$ of $\PP^{v-1}_k$ defined by the semistable points under the action of  $\rho_{\mid SL_{e,k}}$, and the quotient morphism:
$$q: \PP^{v-1}_{k, ss} \lra M(\rho) := \PP^{v-1}_k /\!/ SL_{e,k}.$$
The GIT quotient  $M(\rho)$ is a projective integral $k$-scheme, and is endowed with a canonical  ample  $\Q$-line bundle $L$ such that the following isomorphism of $\Q$-line bundles holds: 
\begin{equation}\label{IntroIsoLonpt}
q^\ast L \simeq 
\cO_{\PP^{v-1}_k}(1)_{\mid \PP^{v-1}_{k, ss}}.
\end{equation}

If $P$ is a $k(C)$-point of $\PP^{v-1}_{k, ss}$, it follows from the valuative criterion of properness that the $k(C)$-point $q_{k(C)}(P)$ of the projective $k$-scheme $M(\rho)$ may be extended into a unique morphism of $k$-schemes:
$$\overline{q_{k(C)}(P)} : C \lra  M(\rho).$$
This allows us to define the \emph{GIT height} of the  point $P$ to be the rational number:
$$\hgt_{GIT}(P) := \deg_C(\overline{q_{k(C)}(P)}^\ast L ).$$

In the arithmetic setting, when the function field $k(C)$ is replaced by a number field, this GIT height is well-documented. In the case where the representation $\rho$ is the natural action of $SL_{N+1,k}$ on effective Chow cycles of fixed dimension and degree in $\PP^N_k$, the idea of using geometric invariant theory to study natural heights of Chow cycles already appears (in both the arithmetic and the geometric setting) in the proof of \cite[Theorems I and III]{Bost94}. The GIT height of Chow cycles is defined in \cite{Zhang96} under the notation $\hat{h}$ and in \cite[3.1]{Bost96} under the notation $h_{\mathrm{inv}}$. This line of study is generalized to more general representations in \cite{Gasbarri00}, and the GIT height is studied for general representations in  \cite[Chapter 4]{Maculan17} under the notation $(1/D) h_{\mathscr{M}_D} \circ \pi$.

\subsubsection{} \label{intro twist}
To a vector bundle $E$ of rank $e$ over $C$ can be attached the twisted vector bundle $E^\rho$ (see for instance \cite{Bogomolov78}, or Subsection \ref{subsection twist} below), and the projective bundle $\PP(E^\rho)$ over $C$. Observe that~$E^\rho$ is a vector bundle of rank $v$, and that the slopes:
$$\mu(E) :=   (\deg_C E)/e
 \quad \mbox{and} \quad \mu(E^\rho) := (\deg_C E^\rho) /v$$
are related by the equality:
$$\mu(E^\rho) = \gamma \, \mu(E).$$

Moreover, for every point $x$ in $C(k)$, the choice of a $k$-basis of $E_x$ determines a basis of $E^\rho_x$, and therefore an isomorphism $\PP(E^\rho)_x \simeq \PP^{v-1}_k.$ The open subscheme $\PP^{v-1}_{k, ss}$ is the image by this isomorphism of an open subscheme $\PP(E^\rho)_{ss,x}$ in  $\PP(E^\rho)_x$, and the morphism $q$ induces, by composition with the restriction of this isomorphism to $\PP(E^\rho)_{ss,x}$, a $GL(E_x)$-invariant morphism:
$$q_x: \PP(E^\rho)_{ss,x} \lra M(\rho).$$
This open subscheme and this morphism do not depend on the choice of a basis of $E_x$. 
As the point $x$ in $C(k)$ varies, these constructions  define a Zariski open subscheme $\PP(E^\rho)_{ss}$ in $\PP(E^\rho)$ and a morphism:
$$q_C: \PP(E^\rho)_{ss} \lra M(\rho).$$ 

Moreover, by using the isomorphism \eqref{IntroIsoLonpt}, one may construct a canonical isomorphism of $\Q$-line bundles over $\PP(E^\rho)_{ss}$:
\begin{equation}\label{IntroIsoLonC}
q_C^\ast L \simeq 
\cO_{E^\rho}(1)_{\mid \PP(E^\rho)_{ss}} \otimes \pi_{\mid \PP(E^\rho)_{ss}}^{\ast}(\det E)^{\otimes \gamma/e},
\end{equation}
where $\pi: \PP(E^\rho) \ra C$ denotes the projection morphism.

\subsubsection{} Let us denote by $\eta$ the generic point of $C$. 

If $V$ is a vector bundle of positive rank $v$ on $C$ whose generic fiber $V_\eta$ is trivialized, every $k(C)$-point~$P$ in $\PP^{v-1}_{k(C)} \simeq \PP(V_\eta)$ can be extended into a unique section
$\overline{P} : C  \ra \PP(V)$ of the projective bundle~$\PP(V)\ra C$. This allows us to define the height $\hgt_{\PP(V)}(P)$ of the point $P$ by the following equality:
\begin{equation}\label{htPV}
\hgt_{\PP(V)}(P) := \deg_C(\overline{P}^\ast \cO_V(1) ) + \mu (V),
\end{equation}
where as above $\mu(V):= \deg V / v$ denotes the slope of the vector bundle $V.$\footnote{The additive normalization by $\mu (V)$ in this definition has the same purpose as the additive normalization in the definition of the height $\hgt_{int}$ in \cite[(1.4.4)]{Mordant22}; see also \eqref{pref hintdef} below. It ensures that the height $\hgt_{\PP(V)}(P)$ is unchanged when the vector bundle $V$ is replaced by $V \otimes L$, where $L$ denotes a line bundle on $C$ whose generic fiber is trivialized: the height $\hgt_{\PP(V)}(P)$ only depends on the projective bundle $\PP:= \PP(V)$, and not on the choice of a vector bundle $V$ such that $\PP \simeq \PP(V)$.}

If $V$ is the twisted vector bundle $E^\rho$ associated to a vector bundle $E$ of rank $e$ as in \ref{intro twist}, a trivialization of the generic fiber $E_\eta$ induces a trivialization of the generic fiber $V_\eta = E^\rho_\eta$. We may therefore attach to a $k(C)$-point $P$ in $\PP^{v-1}_{k(C), ss} \simeq \PP(E^\rho_\eta)_{ss}$ its height $\hgt_{\PP(E^\rho)}(P)$ defined 
by the formula \eqref{htPV} where $V:= E^\rho$, namely the rational number:
$$\hgt_{\PP(E^\rho)}(P) := \deg_C(\overline{P}^\ast \cO_{E^\rho}(1) ) + \gamma \, \mu(E).$$

To a point $P$ in $\PP^{v-1}_{k, ss}(k(C))$ are therefore attached two heights: $\hgt_{GIT}(P)$ and $\hgt_{\PP(E^\rho)}(P)$. 

They are related by the following result, which we will show as a consequence of Theorem \ref{intro height GIT and section} in Subsections \ref{subsection ht rel} and \ref{proof th eq hgt} below.

\begin{theorem} \label{th GITPE} 
For every point $P$ in $\PP^{v-1}_{k, ss}(k(C))$, the following inequality holds:
\begin{equation}\label{GITPE}
\hgt_{GIT}(P) \leq \hgt_{\PP(E^\rho)}(P).
\end{equation}
Moreover equality holds in \eqref{GITPE} if and only the unique extension $\overline{P}$ of $P$ to a section of $\PP(E^\rho)$ over~$C$ 
takes its values in the Zariski open subscheme $\PP(E^\rho)_{ss}$ of $\PP(E^\rho)$.
\end{theorem}

In the arithmetic setting, the comparison of GIT heights and more natural heights analogous to~$\hgt_{\PP(E^\rho)}$ is also well-known. Bost's results \cite[Theorem I and III]{Bost94} in the case of \emph{semistable} Chow cycles essentially state that these natural heights are lower bounded in the arithmetic case and non-negative in the geometric case: the latter property is a straightforward consequence of inequality \eqref{GITPE} and of the obvious non-negativity of $\hgt_{GIT}$. The arithmetic analogue of inequality \eqref{GITPE} is then proved in the case of Chow cycles in \cite[Proposition 4.2]{Zhang96} and in \cite[Proposition 3.1]{Bost96},\footnote{The height associated to Chow cycles in these references is defined differently from our height $\hgt_{\PP(E^\rho)}$; the equality between the two heights (up to multiplication by a constant factor) follows for instance from \cite[Proposition 1.2]{Bost94}.} and is used to show a Northcott property for the natural height in \cite[Proposition 3.2]{Bost96}. For more general representations, this analogue is implicit in \cite{Gasbarri00} and an immediate consequence of \cite[Chapter 4, Theorem 1.5]{Maculan17}.

In the remainder of this paper, we shall mainly be interested in the equality case in Theorem~\ref{th GITPE}, which provides a means of computing the GIT height $\hgt_{GIT}(P)$ provided one has constructed a ``model'' $\overline{P}$ of the point $P$ which is everywhere semistable.  

The fact that \eqref{GITPE} becomes an equality when $\overline{P}$ takes its values in $\PP(E^\rho)_{ss}$ is actually a straightforward consequence of the existence of the ``relative quotient morphism'' $q_C: \PP(E^\rho)_{ss} \ra M(\rho)$ and of the isomorphism \eqref{IntroIsoLonC}.

\subsubsection{}\label{trailerMordantGMS} It turns out that the GIT height $\hgt_{GIT}(P)$ of some $k(C)$-point $P$ of $\PP^{v-1}_{k,ss}$ admits an alternate definition, which is natural in terms of models of $P$ over the curve $C$ and its finite coverings, and does not explicitly involve invariant theory.

Indeed we may introduce the \emph{normalized invariant height} $\hgt_{NI}(P)$ attached to a  point  $P \in \PP^{v-1}(k(C))$ to be the infimum:
\begin{equation}\label{defhtNI}
 \hgt_{NI}(P) := \inf_{C', E'} (\deg \sigma)^{-1} \, \hgt_{\PP(E'^\rho)}({P_{k(C')}})  \in [-\infty, +\infty),
 \end{equation} 
where $C'$ denotes a connected smooth projective curve endowed with a finite morphism
$\sigma : C' \ra C,$
and $E'$ a vector bundle of rank $e$ on $C'$ whose fiber on the generic point $\eta'$ of $C'$ is trivialized, which allows us to identify the point $P_{k(C')} \in \PP^{v-1}(k(C'))$ with a $k(C')$-point of $\PP(E'^\rho)_{\eta'}.$

Observe that, by construction, the height $ \hgt_{NI}(P)$ is invariant under the action of $PGL_e(k(C))$ over $\PP^{v-1}(k(C))$ and is compatible with base change of the base curve $C$
 under finite morphisms. Moreover the inequality \eqref{GITPE} in Theorem \ref{th GITPE} may be rephrased as the following inequality, valid when $P$ belongs to $\PP^{v-1}_{k,ss}(k(C))$:
\begin{equation}\label{preGITNI}
\hgt_{GIT}(P) \leq  \hgt_{NI}(P).
\end{equation}
 
By means of a suitable refinement of the semistable replacement theorem involving finite morphisms of curves, it is possible to show that, when $k$ has characteristic zero, the inequality \eqref{preGITNI} is actually an equality. Indeed we have the following theorem: 
\begin{theorem}\label{theorem GIT NI}
When the base field $k$ has characteristic zero,  for every point $P$ in $\PP^{v-1}_{k,ss}(k(C)),$ the following equality of rational numbers holds:
\begin{equation}\label{GITNI}
\hgt_{GIT}(P) = \hgt_{NI}(P),
\end{equation}
and the infimum in the right-hand side of \eqref{defhtNI} is actually a minimum.

Moreover, for every point $P$ in $\PP^{v-1}(k(C)) \smallsetminus \PP^{v-1}_{k,ss}(k(C)),$ the following equality holds: 
\begin{equation}\label{NIinf}\hgt_{NI}(P) = -\infty.
\end{equation}
\end{theorem}

We will not use this result in the present paper, and we defer its proof to \cite{MordantGMS}, where we shall establish it in the more general framework of good moduli spaces attached to algebraic stacks (see \cite{Alper13, AHLH23}), which encompasses geometric invariant theory. 

In the arithmetic setting, Theorem \ref{theorem GIT NI} is established in the case of Chow cycles in \cite[Proposition 4.2]{Zhang96}. The equivalence between the lower boundedness of the height and the semistability is established for more general representations in \cite[Theorem 1]{Gasbarri00}.

\subsection{The GIT height of a pencil of hypersurfaces $H/C$}

The results of Subsection \ref{GIT height quot}  apply notably in the case of projective hypersurfaces. 

Let us assume that the base field is $k = \C,$ consider $N$ and $d$ two positive integers, denote $e:= N+1$ and $v := \binom{N+d}{N}$ and define the morphism of complex linear groups: 
\begin{equation}\label{rhoclassical}
\rho  : GL_{e,\C} \lra GL_{v, \C}, \quad  g \longmapsto S^d ({}^t\!g^{-1}),
\end{equation}
defined by the $d$-th symmetric power of the dual of the standard  representation of $GL_{e,\C}$. It 
is homogeneous of weight $\gamma = -d$.

The representation $\rho$ is nothing but the classical action of $GL_e(\C)$ by change of coordinates on the vector space $\C[X_0, \dots, X_N]_d$ of homogeneous polynomials of degree $d$ in $N+1$ indeterminates,\footnote{aka $(N+1)$-ary $d$-ics.}  which played a central role in classical invariant theory. 

A $\C(C)$-point $P$ in 
$$\PP^{\binom{N+d}{N} - 1}(\C(C)) \simeq \PP(\C(C)[X_0, \dots, X_N]_d)$$ may be interpreted as a hypersurface $H_\eta$ of degree $d$ in~$\PP^N_{\C(C)}$. 
If the $\C(C)$-point $P$ is in the semistable locus $\PP^{\binom{N+d}{N}-1}_{\C, ss}(\C(C)),$ we shall say that 
\emph{the  hypersurface~$H_\eta$} of degree $d$ in~$\PP^N_{\C(C)}$ defined by $P$ \emph{is semistable}. Moreover we shall call 
the GIT height  of $P$ defined in terms of the representation $\rho$  
the \emph{GIT height of the hypersurface $H_\eta$ in~$\PP^N_{\C(C)}$}, and denote it by:
$$\hgt_{GIT}(H_\eta) := \hgt_{GIT}(P).$$

\subsubsection{} For every vector bundle $E$ of rank $N+1$ on $C,$ the twisted vector bundle $E^\rho$ is canonically isomorphic to $S^d E^\vee,$ and a section $\overline{P}$ of $\PP(E^\rho)$ may be interpreted as a horizontal hypersurface $H$ of relative degree $d$ in $\PP(E)$. 

As in \cite[1.4.2]{Mordant22}, we may attach to such a horizontal hypersurface $H$ its \emph{intersection-theoretic height}, namely the rational number:
\begin{equation}\label{pref hintdef}
\mathrm{ht}_{int}(H/C) := \int_{\PP(E)} c_1(\cO_E(1))^N \cap [H] + d N \mu(E).
\end{equation}

The following result is the consequence of 
 a simple computation of characteristic classes (see the alternate expression of $\hgt_{int}(H/C)$ in terms of a line subbundle of $S^d E^\vee$ defining the hypersurface given in \cite[(6.2.13)]{Mordant22}). 
 
 \begin{proposition}\label{htintsec} Consider $E$ a vector bundle of rank $N+1$ on $C$ whose generic fiber $E_\eta$ is trivialized,  
 $P$ a $\C(C)$-point of $\PP^{\binom{N+d}{N} - 1}_{\C(C)} \simeq \PP(E^\rho_\eta)$, and $\overline{P}$ the section of $\PP(E^\rho)$ over $C$ extending $P$. As above, consider $H$ the horizontal hypersurface of relative degree $d$ in $\PP(E)$ attached to $\overline{P}$.  Then the following equality of rational numbers holds:
$$\hgt_{\PP(E^\rho)}(P) = \hgt_{int}(H/C).$$
\end{proposition}

Using this notation, when applied to the representation \eqref{rhoclassical}, Theorem \ref{th GITPE} translates into the following result.

\begin{theorem}\label{theorem hGITH}
Let $H_\eta$ be a semistable hypersurface of degree $d$ in $\PP^N_{\C(C)} \simeq \PP(E_\eta)$, and let $H$ be the horizontal hypersurface of relative degree $d$  in $\PP(E)$ defined as its closure.  

The following inequality holds:
\begin{equation}\label{hGITH}
\hgt_{GIT}(H_\eta) \leq \hgt_{int}(H/C).
\end{equation}
Moreover equality holds in \eqref{hGITH}  if and only if all the fibers $H_x,$ $x \in C$, seen as hypersurfaces of degree $d$ in the projective spaces $\PP(E_x) \simeq  \PP^N_\C,$ are semistable.
\end{theorem}

\subsubsection{}\label{subsubsection NI for hyp} We may also translate the content of Paragraph \ref{trailerMordantGMS} to the present setting. 

Consider an hypersurface $H_\eta$ of degree $d$ in $\PP^N_{\C(C)}$. We may define its normalized invariant height~$\hgt_{NI}(H_\eta)$ to be the infimum:
\begin{equation}\label{defhtNIbis}
\hgt_{NI}(H_\eta) := \inf_{C', E'} (\deg \sigma)^{-1} \, \hgt_{int}(H'/C') \in [-\infty, +\infty),
\end{equation}
where $C'$ denotes a connected smooth projective curve endowed with a finite morphism
$\sigma : C' \ra C,$ $E'$ a vector bundle of rank $e$ on $C'$ whose fiber on the generic point $\eta'$ of $C'$ is trivialized, and $H'$ the horizontal hypersurface of relative degree $d$ in $\PP(E')$ whose generic fiber $H'_{\eta'}$ coincides with the base change $H_{\eta'}$ of $H_\eta$.

Then  Theorem \ref{theorem GIT NI} becomes the following statement:

\begin{theorem}\label{theorem GIT NIbis}
For every semistable hypersurface $H_\eta$ of degree $d$ in $\PP^N_{\C(C)}$, the following equality of rational numbers holds:
\begin{equation}\label{GITNIbis}
\hgt_{GIT}(H_\eta) = \hgt_{NI}(H_\eta),
\end{equation}
and the infimum in the right-hand side of \eqref{defhtNIbis} is actually a minimum.

Moreover, for every non-semistable hypersurface $H_\eta$ of degree $d$ in $\PP^N_{\C(C)}$, the following equality holds: 
\begin{equation}\label{NIinfbis}\hgt_{NI}(H_\eta) = -\infty.
\end{equation}
\end{theorem}

\section{Comparing the stable Griffiths height and the GIT height of pencils of hypersurfaces I.   Pencils with generic singularities}
\label{sec: GKGIT for double pt}

This section is devoted to a detailed presentation of the results described in Paragraph \ref{1.6} of the introduction, namely to the proof of Theorem  \ref{GK GIT generic}. Let us briefly recall its statement.  

In this section, we consider a connected smooth projective complex curve $C$ with generic point~$\eta$.

If $N$ and $d$ denote two  integers satisfying $N\geq 1$ and $d\geq 2$, then to any smooth hypersurface $H_\eta$ of  degree $d$ in $\PP^N_{\C(C)}$ are attached two natural heights: firstly, the stable Griffiths height~$\hgt_{GK, stab}\big(\H^{N-1}(H_\eta/C_\eta)\big)$ of its middle-dimensional cohomology; and secondly, the GIT height $\hgt_{GIT}(H_\eta)$.

\begin{theorem} \label{theorem GK GIT for double points} Let $N$ and $d$ be two integers satisfying $N \geq 2$ and $d \geq 3$. 

Consider a projective bundle $\pi: \PP(E)\ra C$ defined by a vector bundle $E$ of rank $N+1$ over the smooth projective complex curve $C$, and a hypersurface $H$ in $\PP(E)$ such that $\pi_{\mid H}: H \ra C$ is a flat morphism and the degree
of  the fibers of $\pi_{\mid H}$ is $d$.  

If $H$ defines a pencil of hypersurfaces with generic singularities, namely if $H$ is smooth and if~$\pi_{\mid H}$ admits non-degenerate critical points only, then the following equality holds:
\begin{equation}\label{GKGIT}
\mathrm{ht}_{GK, stab}\big(\H^{N-1}(H_\eta/C_\eta)\big)
= F_{stab} (d,N)\,  \hgt_{GIT}(H_\eta),
\end{equation}
where  $F_{stab}(d,N)$ denotes the rational number in $(1/12) \Z$ defined by \eqref{defFstabodd} and \eqref{defFstabeven}.
\end{theorem}

\subsection{Semistability of singular projective hypersurfaces}
\label{semistab hyp reminder}

In order to compute GIT heights of pencils of hypersurfaces by means of Theorem \ref{theorem hGITH}, it is crucial to have at one's disposal criteria ensuring the semistability of hypersurfaces of a given degree. This question has been investigated in  \cite{Lee08} and \cite{Mordant23}, and we now recall the main results of  these articles.

All these results rely on the Hilbert-Mumford criterion, which was first introduced by Hilbert in \cite{Hilbert93}, then in its modern form by Mumford in \cite[Theorem 2.1]{MumfordFogartyKirwan94}, and has been used by Hilbert, Mumford, and later authors
to investigate the semistability of projective hypersurfaces of given dimension and degree. We refer to \cite{Mordant23} for more precise references.

\subsubsection{}Consider two integers $N \geq 1$ and $d \geq 2$, an algebraically closed field~$k$, and an hypersurface~$H$ of degree $d$ in the projective space  $\PP^N_k$. 

Recall that $H$ is said to be stable (resp. semistable) when the point $[F]$ in $\PP^{\binom{N+d}{N} -1}(k)$ attached to any form $F \in k[X_0, \dots,X_N]_d$ whose vanishing defines $H$ is stable (resp. semistable) relatively to the action of $SL_{N+1, k}$ on $k[X_0, \dots,X_N]_d$ by change of coordinates.

The following theorem is the main result of \cite{Mordant23}. Its proof relies on a remarkable estimate of Benoist \cite[Lemma 3.2]{Benoist14} relating the dimension of the singular locus $H_{\mathrm{sing}}$ of $H$ and the numerical data involved in the Hilbert-Mumford criterion applied to the form $F$. 

\begin{theorem}[{\cite[Theorem 1.1]{Mordant23}}] \label{MainThmGIT}  Let $\delta$ be the maximal multiplicity of $H$ at a point of  $H(k)$, and let $s$ be the dimension of the singular locus  $H_{\mathrm{sing}}$ of $H$.

\begin{enumerate}
\item \label{HM isol mult} If the following condition holds:
\begin{equation} \label{eq HM isol mult}
d \geq \delta \, \min(N+1, s+3)  \quad (\mbox{resp. }  d > \delta \, \min(N+1, s+3)\, ),
\end{equation}
then  $H$ is semistable (resp. stable).

\item  \label{HM isol not cone} Assume $N\geq 2.$ If for every point $P \in H(k)$ where $H$ has multiplicity $\delta$,  
the projective tangent cone $\PP(C_P H)$ in $\PP(T_P \PP^N_k) \simeq \PP^{N-1}_k$ is not the cone\footnote{When $N=2$, this condition has to be interpreted as follows: \emph{the support of  $\PP(C_P H)$ is not a unique point in $\PP(T_P\PP^2_k) \simeq \PP_k^{1}$}.} over some hypersurface in a projective hyperplane of $\PP_k^{N-1}$, and if the following condition holds:
\begin{equation} \label{eq HM isol not cone}
d \geq (\delta-1) \min(N+1, s+3)  \quad (\mbox{resp. }  d > (\delta-1)  \min(N+1, s+3)\, ),
\end{equation} 
then  $H$ is semistable (resp. stable).
\end{enumerate}
\end{theorem}

\subsubsection{} In the remainder of this paper, we shall be especially interested in singular hypersurfaces $H$ with only \emph{semihomogeneous singularities}, namely in hypersurfaces $H$ such that for every  point $P$  in~$H(k)$, the projective tangent cone~$\PP(C_P H)$ is a \emph{smooth} hypersurface in $\PP(T_P \PP^N_k) \simeq \PP^{N-1}_k$. 

Observe that this implies that every  singular point $P$ of $H$ is an isolated singularity (in other words, with the notation of Theorem \ref{MainThmGIT}, $s = 0$), and that when $N \geq 2$, the projective tangent cone~$\PP(C_P H)$ is not the cone over some hypersurface in a projective hyperplane of $\PP_k^{N-1}$. 

\begin{theorem}[{\cite[5.4.2]{Mordant23}, \cite[Proposition 2.5]{Lee08}}]\label{criterionsemihom} Assume that the hypersurface $H$ has only semihomogeneous singularities, and let $\delta$ be the maximum of the multiplicities of the points of $H(k)$. If one of the following two conditions holds:
\begin{equation} \label{cond Mor ss semihom}
N \geq 2 \quad \mbox{and} \quad d \geq 3 (\delta-1),
\end{equation}
or
\begin{equation} \label{cond Lee ss semihom}
d \geq N+1 \quad \mbox{and} \quad  d \geq \delta \, (1 + 1/N),
\end{equation}
then $H$ is semistable.
\end{theorem} 

The semistability of $H$ when \eqref{cond Mor ss semihom} holds follows from the previous observation and Part \eqref{HM isol not cone} of Theorem \ref{MainThmGIT}. Its semistability when \eqref{cond Lee ss semihom} holds follows from Lee's criterion \cite[Proposition 2.5]{Lee08}
involving the log canonical threshold of a singular point; see \cite[Section 5]{Mordant23} for a detailed discussion. 

The interest of criterion \eqref{cond Mor ss semihom}, compared to Lee's criterion \eqref{cond Lee ss semihom},  is to cover the situation when the inequality~$d \leq N$ holds.

In the special case where $\delta$ is at most 2, Theorem \ref{criterionsemihom} becomes:

\begin{corollary}\label{odp} An hypersurface $H$ whose only singularities are ordinary double points is semistable when $N\geq 2$ and $d\geq 3$.
\end{corollary}  

In particular, when $N\geq 2,$ smooth hypersurfaces of degree $d\geq 3$ are semistable. Actually, as shown by Mumford  in \cite[Chapter 4, \S 2]{MumfordFogartyKirwan94}, these are always stable, and smooth quadrics also are semistable. 

As already mentioned, the interested reader may find  a more detailed discussion of (semi)stability of projective hypersurfaces and additional references in \cite{Mordant23}.

\subsection{Pencils of hypersurfaces with  generic singularities}\label{GIT pencil generic}
One of our main results in
 \cite{Mordant22} is the computation of the upper, lower, and stable Griffiths heights of the middle-dimensional cohomology of pencils of projective hypersurfaces, in the ``generic'' situation where the total space is smooth and the singularities of the fibers are  ordinary double points.  
 
 Notably, concerning the stable Griffiths height, we have shown:

\begin{theorem}[{\cite[Theorem 1.4.2]{Mordant22}}]
\label{pref intro GK hyp P(E)}
Let $E$ be a vector bundle of rank~$N+1$ over $C$, and $H \subset \PP(E) $ an horizontal hypersurface of relative degree $d$, smooth over $\C.$  If $\pi_{\mid H}$ has only finitely many critical points, all of which are non-degenerate, then the cardinality of the set $\Sigma$ of critical points satisfies:
\begin{equation}\label{cardSigmahypproj}
\vert \Sigma \vert = (N+1) (d-1)^N \, \mathrm{ht}_{int}(H/C).
\end{equation}

Moreover, under the same hypotheses, the following equality holds:
\begin{equation}\label{HypStab}
\mathrm{ht}_{GK,stab}\big(\H^{N-1}(H_\eta/C_\eta)\big)
= F_{stab}(d,N) \,  \mathrm{ht}_{int}(H/C),
\end{equation}
where 
$F_{stab}(d,N)$ is the element of $(1/12) \Z$ given when  $N$ is odd by:
$$F_{stab}(d,N) := \frac{N+1}{24 d^2} \left[ (d-1)^N  (d^2 N - d^2 - 2 d N - 2 )+ 2 (d^2-1) \right] ,$$
and when $N$ is even by:
$$
F_{stab}(d,N) := \frac{N+1}{24 d^2} \left [ (d-1)^N  (d^2 N + 2 d^2 - 2 d N - 2) - 2 (d^2-1) \right ]. $$
\end{theorem}

Observe that, with the notation of Theorem \ref{pref intro GK hyp P(E)}, all the fibers $H_x$, $x \in C$, are semistable hypersurfaces in $\PP(E_x)$ when $N \geq 2$ and $d\geq 3$, as shown by Corollary \ref{odp}. Consequently we may apply the equality case in Theorem \ref{theorem hGITH} and we obtain that the GIT height of the generic fiber $H_\eta$ satisfies:
\begin{equation}\label{hGITHbis}
\hgt_{GIT}(H_\eta) = \hgt_{int}(H/C).
\end{equation}

By combining the expressions \eqref{HypStab} and \eqref{hGITHbis} for the stable Griffiths and GIT heights, we finally obtain Theorem \ref{theorem GK GIT for double points}.

 \section{Comparing the stable Griffiths height and the GIT height of pencils of hypersurfaces II.  Estimates for  non-generic singularities}
 \label{sec: GKGIT for semi hom}

 \subsection{Pencils of hypersurfaces with semihomogeneous singularities} 
\label{section hyp semihom}

We may want to understand whether the equality \eqref{GKGIT} relating the  stable Griffiths and GIT heights of projective hypersurfaces over $\C(C)$ still holds under more general assumptions concerning the bad reduction of their model over $C$ than the ones in Theorem \ref{theorem GK GIT for double points}. 

To achieve this, we have extended the computations of \cite{Mordant22} concerning stable Griffiths heights, whose outcomes are summarized  in Theorem  \ref{pref intro GK hyp P(E)}, to pencils of hypersurfaces whose singular fibers have at most one semihomogeneous singularity.

Namely, in the second part \cite{MordantGIT2} of this paper, we prove the following theorem.

\begin{theorem}\label{intro GK hyp P(E) hom crit}
Let $C$ be a connected smooth projective complex curve with generic point $\eta$, $E$ a vector bundle of rank~$N+1 \geq 2$ over $C$, and $H \subset \PP(E) $ an horizontal hypersurface of relative degree~$d\geq 2$, smooth over $\C$.
Let us assume that the set of critical points $\Sigma$ of the restriction 
$$f := \pi_{\mid H} : H \lra C$$
is finite, and that the restriction $$\pi_{| \Sigma} : \Sigma  \lra C$$  is injective. 
For every point $P$ in $H$, let $\delta_P$ be the multiplicity at $P$ of the projective hypersurface~$H_{\pi(P)}$ in $\PP(E_{\pi(P)})$.\footnote{Observe that the positive integer $\delta_P$ is at least 2 if and only if $P$ is in $\Sigma.$}
Let us further assume that the singularities of the fibers of $f$ are semihomogeneous, or equivalently that for every point $P$ of $H$ the projective tangent cone $\PP(C_P H_{\pi(P)})$ of~$H_{\pi(P)}$ at $P$ is smooth.

Then the following equality of integers holds:
\begin{equation} \label{intro f(sigma) P(E) hom crit}
\sum_{P \in \Sigma} (\delta_P - 1)^N = (N+1) (d-1)^N \hgt_{int}(H/C).
\end{equation}

Moreover the following equality of rational numbers holds:
\begin{equation} \label{intro GKstabXL P(E) hom crit}
\hgt_{GK, stab}\big(\H^{N-1}(H_\eta/C_\eta)\big) 
= - (N+1) w_{N,d} \; \hgt_{int}(H/C) + \sum_{P \in \Sigma} w_{N,\delta_P},
\end{equation}
where for every positive integer $\delta,$ $w_{N,\delta}$ is the rational number in $(1/12) \Z$ defined by:
\begin{equation}\label{intro wdef}
w_{N,\delta} = (\delta-1)\big[(N \delta + 1) (\delta - 1)^{N-1} + (-1)^N (\delta + 1) \big]/(12 \delta^2).
\end{equation}

\end{theorem}

By combining this theorem with an upper bound for the coefficient $w_{N, \delta}$ and with a study of the cases where equality holds, we shall prove the following proposition in Section \ref{sec: ineq GK semi-hom}.

\begin{proposition}\label{intro upper bound GK}
With the notation of Theorem \ref{intro GK hyp P(E) hom crit}, the following inequality of rational numbers holds:
\begin{equation}
\label{intro ineq GK int hom}
\hgt_{GK, stab}\big(\H^{N-1}(H_\eta/C_\eta)\big) 
\leq F_{stab}(d,N) \, \hgt_{int}(H/C).
\end{equation}
Moreover equality holds in \eqref{intro ineq GK int hom} if and only if one of the following conditions holds:
\begin{enumerate}
\item \label{cond N 1}
$N=1$; 
\item  \label{cond N 2} $N=2$ and $\delta_P = 2$ for every $P \in \Sigma$; 
\item \label{cond N 3} $N=3$ and $\delta_P \in \{2,3\}$ for every $P \in \Sigma$;
\item \label{cond N geq 4} $N\geq 4$ and $\delta_P = 2$ for every $P \in \Sigma$.
\end{enumerate}
\end{proposition}

\subsubsection{}

Observe that, with the notation of Theorem \ref{intro GK hyp P(E) hom crit}, all the fibers $H_x$, $x \in C$, are semistable hypersurfaces in $\PP(E_x)$ when one of the conditions \eqref{cond Mor ss semihom} or \eqref{cond Lee ss semihom} of Theorem \ref{criterionsemihom} holds. Consequently we may apply the equality case in Theorem \ref{theorem hGITH} and we obtain that the GIT height of the generic fiber $H_\eta$ satisfies:
\begin{equation}\label{hGITHter}
\hgt_{GIT}(H_\eta) = \hgt_{int}(H/C).
\end{equation}

Similarly to the derivation of Theorem \ref{theorem GK GIT for double points}, we may combine the results of Theorem \ref{intro GK hyp P(E) hom crit} and Proposition \ref{intro upper bound GK} concerning the heights  
$\hgt_{GK, stab}\big(\H^{N-1}(H_\eta/C_\eta)\big)$ and $ \hgt_{int}(H/C)$ 
with equality \eqref{hGITHter}, which establishes the following result.

\begin{corollary}
\label{pref intro GK hyp P(E) hom crit GIT}
Let $H_\eta$ be a smooth hypersurface of degree $d$ in $\PP^N_{\C(C)}$.

If the hypersurface $H_\eta$ admits as ``model'' an horizontal hypersurface $H \subset \PP(E)$ in the projective bundle $\PP(E)$ associated to a vector bundle $E$ of rank $N+1$ on $C$ whose generic fiber is trivialized, satisfying the hypotheses of Theorem \ref{intro GK hyp P(E) hom crit}, and if, with the notation of that theorem, 
the maximum $\delta$ of the $\delta_P$ for every $P$ in $H$ satisfies the following conditions:
$$N \geq 2 \quad \mbox{and} \quad d \geq 3 (\delta-1),$$
or:
$$ d \geq N+1 \quad \mbox{and} \quad  d \geq \delta \, (1 + 1/N),$$
then the following equality of rational numbers holds:
\begin{equation} \label{pref intro GKstabXL P(E) hom crit GIT}
\hgt_{GK, stab}\big(\H^{N-1}(H_\eta/C_\eta)\big) 
= - (N+1) w_{N,d} \; \hgt_{GIT}(H_\eta) + \sum_{P \in \Sigma} w_{N,\delta_P},
\end{equation}
and the following inequality holds:
\begin{equation}
\label{pref intro ineq GK int hom GIT}
\hgt_{GK, stab}\big(\H^{N-1}(H_\eta/C_\eta)\big) 
\leq F_{stab}(d,N) \, \hgt_{GIT}(H_\eta).
\end{equation}
Moreover equality holds in \eqref{pref intro ineq GK int hom GIT} if and only if one of the conditions \eqref{cond N 1}, \eqref{cond N 2}, \eqref{cond N 3} or \eqref{cond N geq 4} of Proposition~\ref{intro upper bound GK} holds.
\end{corollary}

\subsection{A consequence of the lower semicontinuity of $\hgt_{GK, stab}$}\label{ConsLow}

The following result, concerning pencils of projective hypersurfaces with (GIT-)semistable fibers and with sufficiently large invariant~$\hgt_{int}$, is also useful evidence for the discussion of possible comparison results between the stable Griffiths and GIT heights. 

As above, we denote by $C$ a connected smooth projective complex curve, and by $\eta$ its generic point. We denote  the genus of $C$ by $g$, and we let:
$$(2g-2)^\sim:= \begin{cases} 2g -2 \quad \mbox{ if $g>0$,} \\ -1 \quad \quad \mbox{ if $g=0$.}
\end{cases} 
$$
\begin{proposition}\label{generization} Let $H_\eta$ be a smooth hypersurface of degree $d$ in $\PP^N_{\C(C)}$.

If the hypersurface $H_\eta$ admits as ``model'' an horizontal hypersurface $H \subset \PP(E)$ in the projective bundle $\PP(E)$ associated to a vector bundle $E$ of rank $N+1$ on $C$ whose generic fiber is trivialized, such that the following inequality is satisfied:
\begin{equation}\label{ineqHtIntGeneric}
\hgt_{int}(H/C)  > (2g -2)^\sim + d(\mu_{\max}(E) -\mu(E)),
\end{equation}
 then we have:
\begin{equation}
\label{pref intro ineq GK int hom Bis}
\hgt_{GK, stab}\big(\H^{N-1}(H_\eta/C_\eta)\big) 
\leq F_{stab}(d,N) \, \hgt_{int}(H/C).
\end{equation}

If moreover, for every $x \in C$, the fiber $H_x$ of $H$ over $x$ is a semistable hypersurface in $\PP(E_x),$ then we have:
\begin{equation}
\label{pref intro ineq GK int hom GIT Bis}
\hgt_{GK, stab}\big(\H^{N-1}(H_\eta/C_\eta)\big) 
\leq F_{stab}(d,N) \, \hgt_{GIT}(H_\eta).
\end{equation}

\end{proposition} 

As usual, we denote by $$\mu(E):= \frac{\deg_C E}{\rk E} = \frac{\deg_C E}{N+1}$$ the slope of the vector bundle $E$, and by $\mu_{\max}(E)$ its maximal slope, defined as the maximal value of the slopes of the vector subbundles of positive rank of $E$. 

The validity of \eqref{pref intro ineq GK int hom GIT Bis} under the semistability assumption on the fibers $H_x$ follows from \eqref{pref intro ineq GK int hom Bis} and from the equality case in Theorem \ref{theorem hGITH}.

The proof of \eqref{pref intro ineq GK int hom Bis} will be postponed until Subsection \ref{proof generization GIT} below. It will rely on a deformation argument: by using a lower semicontinuity property of the stable Griffiths height (which will be stated and proved in Section \ref{sec: semicontinuity GK} below), we will reduce to  the case of a pencil of hypersurfaces satisfying the smoothness and non-degeneracy assumptions of Theorem \ref{pref intro GK hyp P(E)}.

\subsection{Questions and remarks}
\label{questions GK GIT}  

\subsubsection{}
The results of Subsections \ref{section hyp semihom} and \ref{ConsLow} lead to the following question, which is a variant of the question concerning the validity of the comparison estimate \eqref{ineq GK GIT semihom} relating $\hgt_{GK, stab}\big(\H^{N-1}(H_{\eta} / C_\eta)\big)$ and~$\hgt_{GIT}(H_\eta)$ discussed in Paragraphs \ref{1.7} and \ref{1.8} of the introduction:

\emph{If $H_\eta$ is a non-singular hypersurface in $\PP^N_{\C(C)}$ of degree $d \geq 2$, does the following inequality hold for every ``model'' $H$ of $H_\eta$ in the projective bundle $\PP(E)$ associated to a vector bundle $E$ of rank~$N+1$ over $C$ whose generic fiber $E_\eta$ is trivialized:}
\begin{equation}
\label{ineq question int} 
\hgt_{GK, stab}\big(\H^{N-1}(H_{\eta} / C_\eta)\big) \leq F_{stab}(d,N) \, \hgt_{int}(H/C) \textrm{ ?}
\end{equation}

Indeed, as established in Theorem \ref{pref intro GK hyp P(E)} and Proposition \ref{intro upper bound GK}, the answer to this question is positive when the singularities of the fibers of the ``model'' $H$ of $H_\eta$ satisfy suitable assumptions. Moreover we have shown in Proposition \ref{generization} that the answer is also positive when the invariant $\hgt_{int}(H/C)$ of the ``model'' $H$ is sufficiently large compared to the instability of $E$.

Let us briefly explain how the comparison estimate \eqref{ineq GK GIT semihom} would follow from a positive answer to the question above.

Consider a non-singular hypersurface $H_\eta$ in $\PP^N_{\C(C)}$ and assume that, for every connected smooth projective curve $C'$ with generic point $\eta'$, for every finite morphism $C' \ra C$, and for every ``model''~$H'$ over $C'$ of the base change~$H_{\eta'}$ of $H_\eta$,  inequality \eqref{ineq question int}  still holds, namely:
$$ \hgt_{GK, stab}\big(\H^{N-1}(H_{\eta'} / C'_{\eta'})\big) \leq F_{stab}(d,N) \, \hgt_{int}(H'/C').$$

The compatibility of the height $ \hgt_{GK, stab}\big(\H^{N-1}(H_{\eta'} / C'_{\eta'})\big)$ with finite extensions of the field~$\C(C)$ and the definition of the height $\hgt_{NI}(H_\eta)$ introduced in \ref{subsubsection NI for hyp} then imply:
\begin{equation}
\label{ineq question NI} 
\hgt_{GK, stab}\big(\H^{N-1}(H_{\eta} / C_\eta)\big) \leq F_{stab}(d,N) \, \hgt_{NI}(H_\eta).
\end{equation}
Moreover, since the smooth hypersurface $H_\eta$ of degree $d \geq 2$ is semistable, according to the equality \eqref{GITNIbis} in Theorem \ref{theorem GIT NIbis}, we have:
\begin{equation}\label{question NI GIT}
\hgt_{NI}(H_\eta) = \hgt_{GIT}(H_\eta).
\end{equation}

Finally \eqref{ineq question NI} and \eqref{question NI GIT} imply the estimate \eqref{ineq GK GIT semihom}, namely:
\begin{equation*}
\label{ineq question GIT} 
\hgt_{GK, stab}\big(\H^{N-1}(H_{\eta} / C_\eta)\big) \leq F_{stab}(d,N) \, \hgt_{GIT}(H_\eta).
\end{equation*}

\subsubsection{}
Among the different cases \eqref{cond N 1}, \eqref{cond N 2}, \eqref{cond N 3}, and \eqref{cond N geq 4} of Proposition \ref{intro upper bound GK} where equality holds in \eqref{intro ineq GK int hom} and \eqref{pref intro ineq GK int hom GIT}, the cases where $\delta_P = 2$  for every $P \in \Sigma$, including cases \eqref{cond N 2} and \eqref{cond N geq 4}, are an obvious consequence of Theorems \ref{pref intro GK hyp P(E)} and \ref{theorem GK GIT for double points}. Moreover the case \eqref{cond N 1} where $N = 1$ follows from the fact that the Griffiths heights of a pencil of $0$-folds vanish, and that the coefficient $F_{stab}(d,1)$ vanishes.

The most interesting case is therefore the case \eqref{cond N 3} where $N = 3$ and where $\delta_P \in \{2, 3\}$ for every~$P \in \Sigma$: it is somewhat surprising that when $N = 3$, semihomogeneous singularities of multiplicity~$3$ behave similarly to ordinary double points when it comes to the comparison of the stable Griffiths height and the GIT height.

\subsubsection{}

We can ask whether an inequality between the stable Griffiths height and the GIT height could hold in the  direction opposite from the one of inequality \eqref{ineq GK GIT semihom} in the semistable case.

Namely, if $H_\eta$ is a non-singular hypersurface in $\PP^N_{\C(C)}$ of degree $d \geq 2$, and if $N \geq 2$ and $(N,d) \neq (3,3)$, does the following inequality of rational numbers hold for every ``model'' $H$ of $H_\eta$ in the projective bundle $\PP(E)$ associated to a vector bundle $E$ of rank $N+1$ over $C$ whose generic fiber is trivialized, such that all the fibers $H_x,$ $x \in C$ of $H$ are \emph{semistable}\footnote{e.g. when the hypotheses on the singularities of Theorem \ref{intro GK hyp P(E) hom crit} and the numerical hypotheses of Corollary \ref{pref intro GK hyp P(E) hom crit GIT} are satisfied.}: 
\begin{equation}\label{other ineq height}
\hgt_{GIT}(H_\eta) = \hgt_{int}(H/C) \leq B(d,N) \, \hgt_{GK, stab}\big(\H^{N-1}(H_\eta/C_\eta)\big)
\end{equation} 
for some non-negative constant $B(d,N)$ ?

Observe that the equality between the two leftmost terms of \eqref{other ineq height} holds because of the semistability hypothesis, and because of the equality case in Theorem \ref{theorem hGITH}.

\subsubsection{}\label{not ineq N 1}

In this paragraph, we explain why an inequality of the form \eqref{other ineq height} cannot hold when $N = 1$ and $d \geq 4$ or when $(N,d) = (3,3)$.

When $N = 1$, for every non-singular hypersurface $H_\eta$ in $\PP^1_{\C(C)}$, the stable Griffiths height $\hgt_{GK, stab}\big(\H^0(H_\eta/C_\eta)\big)$ vanishes.

Also recall that it follows from the Hilbert-Mumford criterion in the case of binary $d$-ics (see \cite[Chapter 4, \S 1]{MumfordFogartyKirwan94}) that a hypersurface in $\PP^1_\C$ of degree $d$ is semistable if and only if it contains no point of multiplicity $> d/2$. 

If moreover $d \geq 4$, then it is possible to find non-singular horizontal hypersurfaces given by cyclic coverings $H \ra C$ with arbitrarily many ramified points, and such that all the ramified indexes are equal to $2$. Let $H$ be a horizontal hypersurface given by such a covering, let $H_\eta$ be its generic fiber, and let $\Sigma$ be its set of critical points, all of which are non-degenerate by hypothesis. 

Using equality \eqref{cardSigmahypproj}, the height 
$$\hgt_{GIT}(H_\eta) = \hgt_{int}(H/C) = \big (2 (d-1)^N \big )^{-1} \vert \Sigma \vert$$ 
will be arbitrarily large while the height $\hgt_{GK, stab}\big(\H^0(H_\eta/C_\eta)\big)$ vanishes, even though the fibers of~$H$ are semistable due to the discussion above. Consequently an inequality of the form \eqref{other ineq height} cannot possibly hold in this case.

In the case where the pair of integers $(N,d)$ is $(3,3),$ using Poincar\'e duality and the fact that the Hodge numbers $h^{0,2}$ and $h^{2,0}$ of a non-singular projective cubic surface vanish, the stable Griffiths height $\hgt_{GK, stab}\big(\H^2(H_\eta/C_\eta)\big)$ of a non-singular cubic hypersurface $H_\eta$ in $\PP^3_{\C(C)}$ always vanishes.

Also recall that it follows from the Hilbert-Mumford criterion in the case of quaternary cubics (see for instance \cite[Chapter 4, \S 2]{MumfordFogartyKirwan94} or \cite[Proposition 6.5]{Beauville09}) that a cubic surface is semistable if and only if it is non-singular or admits only ordinary double points as singularities.

It is possible to find pencils of cubic surfaces on $C$ whose invariant $\hgt_{int}$ is arbitrarily large and whose singular fibers admit only ordinary double points as singularities (and are therefore semistable). This follows for instance from the genericity result \cite[Th\'eor\`eme A.1.3]{MordantMem23}, recalled in Theorem \ref{th genericity in GIT} below, and from the alternate expression of $\hgt_{int}(H/C)$ given in \cite[(6.2.13)]{Mordant22} and recalled in \eqref{htint M E in GIT deformation} below.

Using the semistability property and the equality case in Theorem \ref{theorem hGITH}, the height $\hgt_{GIT}(H_\eta)$ of the generic fiber of such a pencil of cubic surfaces $H$ will be arbitrarily large while the height $\hgt_{GK, stab}\big(\H^2(H_\eta/C_\eta)\big)$ vanishes, even though the fibers of~$H$ are semistable. Consequently an inequality of the form \eqref{other ineq height} cannot possibly hold in this case.

\subsubsection{}Observe finally that there exist values of the pair $(N,d)$ for which  the stable Griffiths height of~$H_\eta$ necessarily vanishes, and  also the GIT height $\hgt_{GIT}(H_\eta)$  necessarily vanishes under the semistability assumption on the fibers of $H$. This is indeed the case when 
$N = 1$ and $d \leq 3$, or  when $d = 2$, as the reader may easily show. Then inequality \eqref{other ineq height} holds trivially.

\section{Relative geometric invariant theory and heights}\label{sec:RelGIT}

This section is devoted to a detailed presentation of the results in Section \ref{GIT hyp} and to their proofs.

Let $k$ be an algebraically closed field, and $e$ and $v$ two positive integers, and let
$$\rho : GL_{e,k} \lra  GL_{v,k}$$
be a morphism of algebraic $k$-groups.

\subsection{The GIT quotient $M(\rho) :=  \PP^{v-1}_k /\!/ SL_{e,k}$}\label{def GIT on pt}

In this section, we review some basic facts concerning the geometric invariant theory associated to the restriction  $\rho_{\mid SL_{e,k}}$ of the representation $\rho$ to the special linear group $SL_{e,k}$.
  
\subsubsection{} \label{def ss us on pt} The $k$-group $SL_{e,k}$ acts on the $k$-scheme $\PP^{v-1}_k$ through the composite morphism 
$$SL_{e,k} 
 \overset{\rho_{\mid SL_{e,k}}}{\lra} GL_{v,k} \lra PGL_{v,k},$$
and the line bundle $\cO(1)$ over $\PP^{v-1}_k$ admits a natural linearization for this action. Consequently~$SL_{e,k}$ acts on the graded algebra of global sections of powers of this line bundle:
$$S^\bullet k^{v \vee} = \bigoplus_{i \geq 0} \Gamma(\PP^{v-1}_k, \cO(i)).$$

Let us consider    its graded sub-algebra of invariant sections:
\begin{equation}\label{InvSubAlg}
\bigoplus_{i \geq 0} \Gamma(\PP^{v-1}_k, \cO(i))^{SL_{e,k}} \subseteq  \bigoplus_{i \geq 0} \Gamma(\PP^{v-1}_k, \cO(i)). 
\end{equation}
It is a finitely generated $k$-algebra (this is a classical result shown by Hilbert in 
\cite{Hilbert90}). 
Therefore, for any sufficiently divisible integer $\delta \geq 1,$ the sections in $\Gamma(\PP^{v-1}_k, \cO(\delta))^{SL_{e,k}}$ generate the $k$-algebra 
$$\bigoplus_{i \geq 0} \Gamma(\PP^{v-1}_k, \cO(i \delta))^{SL_{e,k}}.$$ 

 Let $\delta$ be such an integer and let $\{I_1, \dots, I_p\}$ be a basis  of the $k$-vector space $\Gamma(\PP^{v-1}_k, \cO(\delta))^{SL_{e,k}}$. 

 Let $\PP^{v-1}_{k, us, \delta}$ be the closed subscheme of $\PP^{v-1}_k$ defined as the vanishing locus of the invariant polynomials $I_1, \dots, I_p$. The support $|\PP^{v-1}_{k, us, \delta}|$ is precisely the intersection of the supports of the base loci of the linear systems $\Gamma(\PP^{v-1}_k, \cO(i))^{SL_{e,k}}$ for $i > 0$, and therefore does not depend on the choice of $\delta$. The $k$-points of this support are called the \emph{unstable} $k$-points of   $\PP^{v-1}_k$ endowed with this linearized action of $SL_{e,k}$.

 The \emph{semistable locus} of   $\PP^{v-1}_k$ endowed with this linearized action of $SL_{e,k}$ is defined to be its open subscheme, complement of the closed subscheme $\PP^{v-1}_{k, us, \delta}$:
 $$\PP^{v-1}_{k,ss} := \PP^{v-1}_k \setminus \PP^{v-1}_{k, us, \delta} = \PP^{v-1}_k \setminus (I_1 = \dots = I_p = 0) \subseteq \PP^{v-1}_k.$$
It is also the complement of the intersection of the base loci of the linear systems $\Gamma(\PP^{v-1}_k, \cO(i))^{SL_{e,k}}$ for $i > 0.$

Consequently the inclusion of graded algebras \eqref{InvSubAlg}
defines a $k$-morphism:
$$q : \PP^{v-1}_{k,ss} \lra \PP^{v-1}_k /\!/ SL_{e,k} := \Proj \bigoplus_{i \geq 0} \Gamma(\PP^{v-1}_k, \cO(i))^{SL_{e,k}}.$$
The closed subscheme $\PP^{v-1}_{k, us, \delta}$ and the open subscheme $\PP^{v-1}_{k,ss}$ are stable under the action of $SL_{e,k}$, and the morphism $q$ is $SL_{e,k}$-invariant.

For conciseness' sake, we shall  use the following notation: 
$$M(\rho) := \PP^{v-1}_k /\!/ SL_{e,k}.$$

\subsubsection{}\label{412} The morphism $q$ admits the following more concrete description. 

There is a canonical $k$-isomorphism:
\begin{equation}\label{iso proj mult delta}
 M(\rho)  
 := \Proj \bigoplus_{i \geq 0} \Gamma(\PP^{v-1}_k, \cO(i))^{SL_{e,k}} \simeq \Proj \bigoplus_{i \geq 0} \Gamma(\PP^{v-1}_k, \cO(i \delta))^{SL_{e,k}},
 \end{equation}
and the canonical surjective morphism of graded $k$-algebras:
$$S^\bullet \Gamma(\PP^{v-1}_k, \cO(\delta))^{SL_{e,k}} \lra \bigoplus_{i \geq 0} \Gamma(\PP^{v-1}_k, \cO(i \delta))^{SL_{e,k}}$$
induces a closed immersion of $k$-schemes:
\begin{equation}\label{imm proj sl} \Proj \bigoplus_{i \geq 0} \Gamma(\PP^{v-1}_k, \cO(i \delta))^{SL_{e,k}} \hlra \PP\big(\Gamma(\PP^{v-1}_k, \cO(\delta))^{SL_{e,k} \vee}\big) \lrasim 
 \PP^{p-1}_k,
\end{equation}
where the last isomorphism is induced by the basis $(I_1, \dots,I_p)$ of $\Gamma(\PP^{v-1}_k, \cO(\delta))^{SL_{e,k}}$.

Composed with the isomorphism \eqref{iso proj mult delta} and the closed immersion \eqref{imm proj sl}, the morphism $q$ becomes the morphism defined by the invariant polynomials $I_1,\dots,I_p$:
$$[I_1 :\cdots : I_p] : \PP^{v-1}_{k,ss} \overset{q}{\lra} M(\rho) =\PP^{v-1}_k /\!/ SL_{e,k} \lrasim \Proj \bigoplus_{i \geq 0} \Gamma(\PP^{v-1}_k, \cO(i \delta))^{SL_{e,k}} \,\hlra \PP^{p-1}_k$$
defined by:
$$[I_1 :\cdots : I_p]([x_1 : \cdots: x_v]) :=    [I_1(x_1, \dots, x_v):\dots: I_p(x_1, \dots,x_v)].$$

\subsubsection{}Observe that by construction, there is a canonical isomorphism of line bundles over $\PP^{v-1}_{k,ss}$:
$$q^* \cO_{\PP^{p-1}_k}(1)_{| M(\rho)} \simeq [I_1:\cdots :I_p]^* \cO_{\PP^{p-1}_k}(1) \simeq \cO_{\PP^{v-1}_{k,ss}}(\delta).$$
Consequently, denoting by $L$ the $\Q$-line bundle $\cO_{\PP^{p-1}_k}(1)_{| M(\rho)}^{\otimes 1/\delta}$ on $M(\rho)$, there is an isomorphism of $\Q$-line bundles over $\PP^{v-1}_{k,ss}$:
$$q^* L \lrasim \cO_{\PP^{v-1}_{k,ss}}(1).$$

The isomorphism class of the $\Q$-line bundle $L$ does not depend on the choice of the integer $\delta$ or of the invariant polynomials~$I_1,\dots,I_p.$

\subsection{Homogeneous representations of linear groups, twisting of vector bundles, and geometric invariant theory over a general base}\label{subsection twist}

\subsubsection{}
Recall that, if $S$ is a $k$-scheme, the theory of principal bundles (see for instance \cite{Bogomolov78}) allows one to associate a vector bundle $E^\rho$  of rank $v$ over $S$ to any vector bundle $E$ of rank $e$ over~$S$. 

If $E$ may be trivialized over a Zariski open covering $(U_i)_{i \in I}$ of $S,$ and therefore may be defined by some cocycle $(h_{ij})_{(i,j) \in I^2},$ with 
$$h_{ij}: U_i \cap U_j \lra GL_{e,k}$$
a $k$-morphism, then the vector bundle $E^\rho$ may be defined by the cocycle $(\rho \circ h_{ij})_{(i,j) \in I^2}$ with values in $GL_{v,k}$. This construction is compatible with isomorphisms of vector bundles, and with arbitrary base change of the $k$-scheme $S.$

\subsubsection{} From now on, we  assume that \emph{the morphism of $k$-linear groups $\rho$ is a homogeneous representation}. 
 In other words, we assume that there exists $\gamma \in \Z$  --- the \emph{weight} of $\rho$ --- such that, for every $t$ in $k^\ast$, the following equality holds  in $GL_v(k)$:
$$\rho(t\, \Id_e) = t^{\gamma} \, \Id_v.$$
This implies that the subscheme $\PP^{v-1}_{k,ss}$ of $\PP^{v-1}_{k}$ is stable under the action of $GL_{e,k}$ and that the morphism $q$ is $GL_{e,k}$-invariant. 

It also implies that, for any positive integer $\delta$ such that $\Gamma(\PP^{v-1}_k, \cO(\delta))^{SL_{e,k}} \neq \{0\},$ the following divisibility holds:
$$ e \mid \gamma \delta.$$

The construction of $\PP^{v-1}_{k,ss}$ and of $q$ admits the following relative version valid over an arbitrary $k$-scheme. 

\begin{proposition}\label{prop:relGIT}
To any $k$-scheme $S$ and to any vector bundle $E$ of rank $e$ over $S$,  it is possible to associate, in a unique way, an open subscheme $\PP(E^\rho)_{ss}$ of 
the $S$-scheme $\PP(E^\rho)$ and a $k$-morphism
$$q_{S,E}: \PP(E^\rho)_{ss} \lra M(\rho)$$
such that the following conditions are satisfied:
\begin{enumerate}[(i)]
\item If $S = \Spec k$ and $E = k^{\oplus e}$, then $\PP(E^\rho)_{ss} = \PP^{v-1}_{k, ss}$ and $q_{S,E} = q$.
\item Let $\sigma: S' \ra S$  be a morphism of $k$-schemes, $E$ (resp. $E'$) a vector bundle of rank $e$ over $S$ (resp. over $S'$), and 
$\phi : E' \lrasim \sigma^\ast E$ 
an isomorphism of vector bundles over $S'$. If we denote by 
$$\phi^\rho: E'^\rho \lrasim (\sigma^\ast E)^\rho \simeq \sigma^\ast (E^\rho)$$
and
$$\PP(\phi^\rho) : \PP(E'^\rho) \lrasim \PP(E^\rho) \times_{S, \sigma} S'$$
the associated isomorphisms of vector bundles and of projective bundles over $S'$, 
then the following equality of open subschemes of $\PP(E^\rho) \times_{S, \sigma} S'$ holds:
$$\PP(\phi^\rho) (\PP(E'^\rho)_{ss}) = \PP(E^\rho)_{ss} \times_{S, \sigma} S',$$
and the following diagram is commutative: 
$$
\xymatrix{
\PP(E'^\rho)_{ss} \ar[r]^-{\PP(\phi^\rho)}_-{\sim} 
\ar[dr]_{q_{S',E'}} & \PP(E^\rho)_{ss} \times_{S, \sigma} S' \ar[r] & \PP(E^\rho)_{ss} \ar[dl]^{q_{S,E}} \\
& M(\rho) &
}
$$
where the second horizontal arrow is the first projection. 
\end{enumerate}

Moreover, for every $k$-scheme $S$ and every vector bundle $E$ of rank $e$ over $S$, we have an isomorphism of $\Q$-line bundles over $\PP(E^\rho)_{ss}$, where $\pi: \PP(E^\rho) \ra S$ denotes the structural morphism:
$$q_{S,E}^\ast L \lrasim \big(\cO_{E^\rho}(1) \otimes \pi^\ast (\det E)^{\otimes \gamma/e} \big)_{| \PP(E^\rho)_{ss} }.$$
This isomorphism is canonical up to some roots of unity, and equivariant with respect to the natural $GL(E)$-linearizations of its source and its range.
\end{proposition}

This is a straightforward consequence of  the $GL_{e,k}$-invariance of the morphism $q: \PP^{v-1}_{k, ss} \ra M(\rho)$ and of the definition of the twisted vector bundles $E^\rho$, and we leave the details to the interested reader. Actually, with the notation of \ref{412}, for every $S$ and $E$ as above, there is a canonical isomorphism of line bundles on $\PP(E^\rho)_{ss}$:
\begin{equation}\label{iso line quot delta rel}
q_{S,E}^\ast L^{\otimes \delta} := q_{S, E}^\ast \cO_{\PP^{p-1}_k} (1)_{\mid M(\rho)} \lrasim (\cO_{E^\rho}(\delta) \otimes \pi^\ast (\det E)^{\otimes \gamma \delta /e} )_{| \PP(E^\rho)_{ss} }.
\end{equation}

 This construction already appears, at least implicitly, in the work of Bogomolov \cite{Bogomolov78}. It is a geometric analogue of the construction of the semistable locus and quotient of a projective bundle over the ring of integers of a number field in \cite[Chapter 4, Proposition 3.5]{Maculan17}.

 \subsubsection{}\label{interpretation q rel with delta} The  construction in Proposition \ref{prop:relGIT} may also be interpreted as follows. We adopt the notation of Subsection \ref{def GIT on pt}; in particular, we consider a positive integer $\delta$ such that the sections in $\Gamma(\PP^{v-1}_k, \cO(\delta))^{SL_{e,k}}$ generate the $k$-algebra 
$$\bigoplus_{i \geq 0} \Gamma(\PP^{v-1}_k, \cO(i \delta))^{SL_{e,k}}.$$ 

 Let $S$ be a non-empty $k$-scheme, which for simplicity will be assumed to be reduced, $E$ a vector bundle of rank $e$ on $S$, and $\pi: \PP(E^\rho)\ra S$ the associated projective bundle. Since $\rho$ is homogeneous of weight $\gamma$, we may attach to any $SL_{e,k}$-invariant section $F$ in $\Gamma(\PP^{v-1}_k, \cO(\delta))$ a unique section~$F_{S,E}$ of the line bundle over $\PP(E^\rho)$:
 $$\cO_{E^\rho}(\delta) \otimes \pi^\ast (\det E)^{\otimes \gamma \delta / e},$$
 such that for every point $x$ in $S(k)$ and for every isomorphism of vector spaces $k^{\oplus e} \simeq E_x$, the induced isomorphism of $k$-schemes $\PP^{v-1}_k \simeq \PP(E^\rho_x)$ identifies the restriction $F_{S,E | \PP(E^\rho_x)}$ with $F$. This defines an injective morphism of $k$-vector spaces:
 $$\Gamma(\PP^{v-1}_k, \cO(\delta))^{SL_{e,k}} \lra \Gamma \Big (\PP(E^\rho), \cO_{E^\rho}(\delta) \otimes \pi^\ast (\det E)^{\otimes \gamma \delta / e} \Big).$$
We denote by $W_{S,E, \delta}$ its image. 
 
 The closed subscheme $\PP(E^\rho)_{us, \delta}$ of $\PP(E^\rho)$ defined as the base locus of the linear system $W_{S, E, \delta}$ is a relative analogue of the closed subscheme $\PP^{v-1}_{k, us, \delta}$ introduced in \ref{def ss us on pt}, in the sense that for every point $x$ in $S(k)$ and for every isomorphism of vector spaces $k^{\oplus e} \simeq E_x$, the induced isomorphism of~$k$-schemes $\PP^{v-1}_k \simeq \PP(E^\rho_x)$ identifies the subscheme
 $$\PP(E^\rho_x)_{us, \delta} := \PP(E^\rho)_{us, \delta} \cap \PP(E^\rho_x)$$ of $\PP(E^\rho_x)$ with $\PP^{v-1}_{k, us, \delta}$.

 The support $|\PP(E^\rho)_{us, \delta}|$ does not depend on the choice of $\delta$, and its complement in $\PP(E^\rho)$ is precisely the open subscheme $\PP(E^\rho)_{ss}$.

 Moreover, as in \ref{412}, the composition of the relative quotient morphism $q_{S,E}$ with the closed immersion
 \begin{equation*} i_\delta : M(\rho) \simeq \Proj \bigoplus_{i \geq 0} \Gamma(\PP^{v-1}_k, \cO(i \delta))^{SL_{e,k}} \hlra \PP\big(\Gamma(\PP^{v-1}_k, \cO(\delta))^{SL_{e,k} \vee}\big) \lrasim
 \PP_k(W_{S, E, \delta}^\vee)
\end{equation*}
is precisely the morphism:
$$h_{S,E, \delta} : \PP(E^\rho)_{ss} = \PP(E^\rho) \setminus \PP(E^\rho)_{us, \delta} \lra \PP_k(W_{S, E, \delta}^\vee)$$
induced by the linear system $W_{S, E, \delta}$. 

By definition of $h_{S,E, \delta}$, there is a canonical isomorphism of line bundles on $\PP(E^\rho)_{ss}$:
\begin{equation} \label{first iso GIT rel inv sections}
    h_{S,E, \delta}^\ast \cO_{W_{S, E, \delta}^\vee}(1) \simeq \big(\cO_{E^\rho}(\delta) \otimes \pi^\ast (\det E)^{\otimes \gamma \delta / e}\big)_{| \PP(E^\rho)_{ss} },
\end{equation}
and by definition of $L$, there is a canonical isomorphism of line bundles on $M(\rho)$:
\begin{equation} \label{second iso GIT rel inv sections}
    i_\delta^\ast \cO_{W_{S, E, \delta}^\vee}(1) \simeq L^{\otimes \delta}.
\end{equation}

The isomorphisms \eqref{first iso GIT rel inv sections} and \eqref{second iso GIT rel inv sections}, along with the equality of morphisms
\begin{equation}\label{equality comp immersion quot rel}
    h_{S, E, \delta} = i_\delta \circ q_{S,E} : 
    \PP(E^\rho)_{ss} \lra \PP_k(W_{S, E, \delta}^\vee),
\end{equation}
observed above, yield a different proof of the isomorphism \eqref{iso line quot delta rel}.
 
\subsection{Relative geometric invariant theory over projective  curves and heights}\label{subsection ht rel}

In this subsection, we consider a connected smooth projective curve $C$ with generic point $\eta$ over $k$, and we apply the previous construction when $S$ is $C$ to 
investigate the GIT height of $k(C)$-points of $\PP^{v-1}_{k, ss}$ introduced in Subsection \ref{GIT height quot}. 

Recall that if $P$ is a point in $\PP^{v-1}_{k,ss}(k(C))$, the image $q_{k(C)}(P)$ in $M(\rho)(k(C))$ may be extended into a unique morphism of $k$-schemes:
$$\overline{q_{k(C)}(P)} : C \lra  M(\rho),$$
and that the \emph{GIT height} of the  point $P$ is defined to be the rational number:
$$\hgt_{GIT}(P) := \deg_C(\overline{q_{k(C)}(P)}^\ast L ).$$

Let $E$ be  a vector bundle of rank $e$ on $C$ whose generic fiber $E_\eta$ is trivialized. As in Subsection \ref{GIT height quot}, a point $P$ in $\PP^{v-1}(k(C)) \simeq \PP(E^\rho_\eta)(k(C))$ may be extended into a unique section $\overline{P}$ of the projective bundle $\pi:\PP(E^\rho) \ra C$, and the height $\hgt_{\PP(E^\rho)}(P)$ is defined to be the rational number:
$$\hgt_{\PP(E^\rho)}(P) := \deg_C(\overline{P}^\ast \cO_{E^\rho}(1) ) + \gamma \, \mu(E),$$
where $\mu(E) := (\deg_C E)/e$ denotes the slope of $E$.

As in Subsection \ref{def GIT on pt} and \ref{interpretation q rel with delta}, we consider a positive integer $\delta$  such that the sections in the vector space $\Gamma(\PP^{v-1}_k, \cO(\delta))^{SL_{e,k}}$ generate the $k$-algebra 
$$\bigoplus_{i \geq 0} \Gamma(\PP^{v-1}_k, \cO(i \delta))^{SL_{e,k}}.$$ 

Using a blow-up argument as in \cite[4.4]{Fulton98} and \cite[Lemma 2.3]{Bost96},  we shall compare the two heights $\hgt_{GIT}(P)$ and $\hgt_{\PP(E^\rho)}(P)$ attached to a point $P$ in $\PP^{v-1}_{k, ss}(k(C))$.

Namely we shall show the following result, which admits Theorem \ref{th GITPE} as a consequence.

\begin{theorem}\label{intro height GIT and section}
 With the above notation, for  every point $P$ in $\PP^{v-1}_{k, ss}(k(C))$, the following equality of rational numbers is satisfied:
 \begin{equation}\label{eq:intro height GIT and section}
\hgt_{\PP(E^\rho)}(P) = \hgt_{GIT}(P) + \delta^{-1} \mathrm{lg} \big (\overline{P}^{-1}(\PP(E^\rho)_{us, \delta}) \big),
\end{equation}
where $\mathrm{lg}$ denotes the length of a  closed subscheme of $C$ with finite support, where $\PP(E^\rho)_{us, \delta}$ denotes the closed subscheme of $\PP(E^\rho)$ defined in \ref{interpretation q rel with delta}, and where $\overline{P}$ denotes the section over~$C$ of $\PP(E^\rho)$ extending~$P$. 
\end{theorem}

Observe that the closed subscheme $\overline{P}^{-1}(\PP(E^\rho)_{us, \delta})$ of $C$ is indeed strict, hence with finite support, because the generic fiber~$P$ of $\overline{P}$ lies in $\PP(E^\rho)_{ss, \eta}(k(C))$. 

The second term in the right-hand side of \eqref{eq:intro height GIT and section} is a sum of local terms attached to the points of unstable reduction of the point $P$. These local terms play a very similar role in this equality as the ``instability measures'' introduced in 
\cite[Chapter 4, 1.2]{Maculan17} do in the arithmetic analogue of this theorem (\cite[Chapter 4, Theorem 1.5]{Maculan17}).

Before giving the proof of Theorem \ref{intro height GIT and section} in Subsection \ref{proof th eq hgt}, let us explain how to derive Theorem~\ref{th GITPE} from it.

The second term in the right-hand side of \eqref{eq:intro height GIT and section} is always non-negative, which establishes inequality~\eqref{GITPE}. 
Moreover this second term vanishes if and only if the subscheme $\overline{P}^{-1}(\PP(E^\rho)_{us, \delta})$ of~$C$ is empty, or equivalently, if and only if the section $\overline{P}$ has values in the complement subscheme $\PP(E^\rho) \smallsetminus \PP(E^\rho)_{us, \delta}$. Since this complement is precisely the open subscheme $\PP(E^\rho)_{ss}$, 
this yields the description of  the equality case in inequality \eqref{GITPE} stated in Theorem \ref{th GITPE}.

\subsection{Heights of sections and orders of contact of sections with the subscheme $\PP(E^\rho)_{us, \delta}$} \label{proof th eq hgt}

In this subsection, we establish Theorem \ref{intro height GIT and section}.

\subsubsection{Construction of a blow-up of the projective bundle $\PP(E^\rho)$}

We adopt the notation introduced in~\ref{interpretation q rel with delta}.

Consider the blow-up
$$\nu : Y \lra \PP(E^\rho)$$ of the closed subscheme $\PP(E^\rho)_{us, \delta}$ of $\PP(E^\rho)$, and let $D$ be the Cartier divisor in $Y$ defined as its exceptional divisor.

Observe that $\nu$ induces an isomorphism between $Y \smallsetminus D$ and $\PP(E^\rho)_{ss}$. Consequently the $k(C)$-point $P$ in $\PP_{k,ss}^{v-1}(k(C)) \simeq \PP(E^\rho_\eta)_{ss}(k(C))$ admits a unique preimage $\nu_{\mid Y_\eta \smallsetminus D_\eta}^{-1}(P)$ in $(Y_\eta \smallsetminus D_\eta)(k(C))$. Since the $k$-scheme $Y$ is projective, this $k(C)$-point may be extended into a morphism:
$$\overline{\nu_{\mid  Y_\eta \setminus D_\eta}^{-1}(P)} : C \lra  Y.$$

By unicity of the extension from $\eta$ to $C$ of morphisms to the separated scheme $\PP(E^\rho)$, the following equality of morphisms from $C$ to $\PP(E^\rho)$ holds:
\begin{equation}\label{eq unicity nu for ht GIT sections}
 \nu \circ \overline{\nu_{\mid Y_\eta \setminus D_\eta}^{-1}(P)} = \overline{P}.
 \end{equation}

A morphism of $k$-schemes from $Y$ to $\PP(W_{C, E, \delta}^\vee)$ may be constructed as follows. Let $\cI$ be the ideal sheaf in $\cO_{\PP(E^\rho)}$ associated to the closed subscheme $\PP(E^\rho)_{us, \delta}$. As this closed subscheme is precisely the base locus of the linear system $W_{C, E, \delta}$ of sections of the line bundle $\cO_{E^\rho}(\delta) \otimes \pi^\ast (\det E)^{\otimes \gamma \delta / e}$, there is a canonical surjective morphism of $\cO_{\PP(E^\rho)}$-modules:
\begin{equation}\label{morphismWcI}
W_{C, E, \delta} \otimes_{k} \cO_{\PP(E^\rho)} \lra \cO_{E^\rho}(\delta) \otimes \pi^\ast (\det E)^{\otimes \gamma \delta / e} \otimes \cI.
\end{equation}

Tensoring both sides of \eqref{morphismWcI} by the line bundle $\cO_{E^\rho}(- \delta) \otimes \pi^\ast (\det E)^{\otimes - \gamma \delta / e},$ then taking the graded algebras of symmetric powers and ideal powers, we obtain a surjective morphism of graded~$\cO_{\PP(E^\rho)}$-algebras:
$$S^\bullet \big(W_{C, E, \delta} \otimes_{k} \cO_{E^\rho}(- \delta) \otimes \pi^\ast (\det E)^{\otimes - \gamma \delta / e}\big) \lra \bigoplus_{n \geq 0} \cI^n.$$

This morphism induces through the Proj construction a closed immersion of $\PP(E^\rho)$-schemes:
\begin{align*}Y := \mathrm{Proj}_{\PP(E^\rho)} \bigoplus_{n \geq 0} \cI^n  \,\hlra\, & \mathrm{Proj}_{\PP(E^\rho)} S^\bullet \big(W_{C, E, \delta} \otimes_{k } \cO_{E^\rho}(- \delta) \otimes \pi^\ast (\det E)^{\otimes - \gamma \delta / e}\big) \\
& \quad \quad \simeq \PP_{\PP(E^\rho)} \big(W_{C, E, \delta}^\vee \otimes_{k } \cO_{E^\rho}(\delta) \otimes \pi^\ast (\det E)^{\otimes \gamma \delta / e}\big)\\ 
& \quad \quad \simeq \PP(E^\rho) \times_k \PP(W_{C, E, \delta}^\vee).
\end{align*}
The composition of this closed immersion with the second projection will be denoted by:
$$f_\delta : Y \lra  \PP(W_{C, E, \delta}^\vee).$$

According to the functoriality properties of the Proj  construction,  there is a natural isomorphism of line bundles on $Y$:
\begin{equation}\label{iso f for ht GIT sections}
f_\delta^\ast \cO_{W_{C, E, \delta}^\vee}(1) \simeq \nu^\ast \cO_{E^\rho}(\delta) \otimes \nu^\ast \pi^\ast (\det E)^{\otimes \gamma \delta/e} \otimes \cO_Y(-D).
\end{equation}

Moreover it follows from the  construction of $f_\delta$ and from equality \eqref{equality comp immersion quot rel} that the following equality of morphisms holds:
\begin{equation}\label{eq f for ht GIT sections}
f_{\delta \mid Y \setminus D} = i_\delta \circ q_{C,E} \circ \nu_{| Y \setminus D} : Y \setminus D \lra \PP(W_{C,E,\delta}^\vee).
\end{equation}
Indeed both sides of \eqref{eq f for ht GIT sections} map any $k$-point $x$ in $(Y\smallsetminus D)(k)$ to the hyperplane in the linear system~$W_{C, E, \delta}$ consisting of sections vanishing at the point $\nu(x)$ of $\PP(E^\rho)_{ss}(k).$

Furthermore, using the unicity of the  extension from $\eta$ to $C$ of morphisms to $\PP(W_{C,E,\delta}^\vee)$ and equality \eqref{eq f for ht GIT sections}, we obtain  the following equality of morphisms:
\begin{equation}\label{eq unicity f nu for ht GIT sections} 
f_\delta \circ \overline{\nu_{\mid Y_\eta \setminus D_\eta}^{-1}(P)} = i_\delta \circ \overline{q_{k(C)}(P)}: C \lra \PP(W_{C,E,\delta}^\vee).
\end{equation}

These constructions are summarized in the following diagram, where  the square with horizontal and vertical arrows is cartesian\footnote{The reader should beware that this diagram is not commutative. Compositions of arrows where $C$ occurs not only as the starting or the end object should be excluded from commutativity conditions.}:
\begin{equation*}\label{diagram ht GIT blow up}
\xymatrix{
Y \ar[d]^{\nu} \ar@/^2.5pc/[rrrd]^{f_\delta} & \,Y \setminus D \ar@{_{(}->}[l] \ar[d]_{\rotatebox{90}{$\sim$}}^{\nu_{| Y \setminus D}} & & \\
\PP(E^\rho) \ar[d]^{\pi} & \,\PP(E^\rho)_{ss} \ar@{_{(}->}[l] \ar[r]^{q_{C,E}} 
& M(\rho)\, \ar@{^{(}->}[r]^{i_\delta\;\;} & \PP(W_{C, E, \delta}^\vee) \\
C \ar@/^1pc/[u]^{\overline{P}} \ar@/^3pc/[uu]^{\overline{\nu_{| Y_\eta \setminus D_\eta}^{-1}(P)}} \ar[rru]_{\overline{q_{k(C)}(P)}} & 
& & \quad \quad \quad \quad \quad .
}
\end{equation*}

\subsubsection{Comparison of  line bundles pulled back on $C$}
Using successively the isomorphism \eqref{second iso GIT rel inv sections}, equality \eqref{eq unicity f nu for ht GIT sections}, the isomorphism \eqref{iso f for ht GIT sections}, equality \eqref{eq unicity nu for ht GIT sections}, and the fact that $\overline{P}$ is a section of $\pi$, we obtain the following isomorphisms of line bundles on $C$:
\begin{align*}
\overline{q_{k(C)}(P)}^\ast L^{\otimes \delta} &\simeq \big(i_\delta \circ \overline{q_{k(C)}(P)} \big)^\ast \cO_{W_{C, E, \delta}^\vee}(1) \\
& \simeq  \big(f_\delta \circ \overline{\nu_{\mid Y_\eta \setminus D_\eta}^{-1}(P)} \big)^\ast \cO_{W_{C, E, \delta}^\vee}(1) \\
& \simeq \overline{\nu_{\mid Y_\eta \setminus D_\eta}^{-1}(P)}^\ast \big(\nu^\ast \cO_{E^\rho}(\delta) \otimes \nu^\ast \pi^\ast (\det E)^{\otimes \gamma \delta/e} \otimes \cO_Y(-D) \big) \\
& \simeq \overline{P}^\ast \cO_{E^\rho}(\delta) \otimes \overline{P}^\ast \pi^\ast (\det E)^{\otimes  \gamma \delta/e} \otimes \overline{\nu_{\mid  Y_\eta \setminus D_\eta}^{-1}(P)}^\ast \cO_Y(-D) \\
& \simeq \overline{P}^\ast \cO_{E^\rho}(\delta) \otimes (\det E)^{\otimes \gamma \delta/e} \otimes \overline{\nu_{\mid Y_\eta \setminus D_\eta}^{-1}(P)}^\ast \cO_Y(-D).
\end{align*}
Taking degrees, we obtain the following equality of integers:
$$\delta \, \deg_C \big(\overline{q_{k(C)}(P)}^\ast L\big) = \delta \, \deg_C \big(\overline{P}^\ast \cO_{E^\rho}(1) \big) + \gamma \delta \, \mu(E) + \deg_C\big(\overline{\nu_{\mid Y_\eta \setminus D_\eta}^{-1}(P)}^\ast \cO_Y(-D) \big),$$
which we may rewrite as follows using the definition of the heights $\hgt_{GIT}(P)$ and $\hgt_{\PP(E^\rho)}(P)$:
$$\delta \, \hgt_{GIT}(P) = \delta\, \hgt_{\PP(E^\rho)}(P) + \deg_C \big(\overline{\nu_{\mid Y_\eta \setminus D_\eta}^{-1}(P)}^\ast \cO_Y(-D) \big),$$
or equivalently:
\begin{equation}\label{first eq comparison hgt} 
\hgt_{\PP(E^\rho)}(P) = \hgt_{GIT}(P) + \delta^{-1} \deg_C\big(\overline{\nu_{\mid Y_\eta \setminus D_\eta}^{-1}(P)}^\ast \cO_Y(D) \big).
\end{equation}

Observe that the integer $\deg_C\big(\overline{\nu_{\mid Y_\eta \smallsetminus D_\eta}^{-1}(P)}^\ast \cO_Y(D) \big)$ 
is precisely the length of the strict subscheme $\overline{\nu_{\mid Y_\eta \smallsetminus D_\eta}^{-1}(P)}^{-1}(D)$ in $C.$ By definition of the exceptional divisor $D$ in $Y,$ this subscheme may also be written as follows:
\begin{align}
\overline{\nu_{\mid  Y_\eta \setminus D_\eta}^{-1}(P)}^{-1} (D) & = \big(\nu \circ \overline{\nu_{\mid Y_\eta \setminus D_\eta}^{-1}(P)}\big)^{-1} (\PP(E^\rho)_{us, \delta}) \nonumber \\
\label{eq comparison subschemes} & = \overline{P}^{-1} (\PP(E^\rho)_{us, \delta}),
\end{align}
where in \eqref{eq comparison subschemes} we have used equality \eqref{eq unicity nu for ht GIT sections}.

Replacing in \eqref{first eq comparison hgt} yields the desired equality \eqref{eq:intro height GIT and section} and concludes the proof.

\section{Expression as an integral and lower semicontinuity  of the stable Griffiths height}
\label{sec: semicontinuity GK}

In this section we recall Peters' expression of the Griffiths height of a polarized variation of Hodge structures with unipotent monodromy as an integral of the Chern form on the Griffiths line bundle defined by the Hodge metric. This expression also computes the stable Griffiths height when the monodromy is not assumed to be unipotent. 

Using this result, we will show a lower semicontinuity property for the stable Griffiths height under deformation of the variation of Hodge structures. This property will be used in Subsection \ref{proof generization GIT} below as part of the proof of Proposition \ref{generization}.

\subsection{A theorem of Peters} Let us first recall Peters' result. As above, let $C$ be a connected smooth projective complex curve with generic point $\eta$, let $\Cc$ be the complement in $C$ of a finite subset $\Delta$, and let $\V = (V_\Z, \cF^\bullet)$ be an integral variation of Hodge structures of weight $n$ on $\Cc$, where~$\cF^\bullet$ denotes the Hodge filtration on the vector bundle
$$\cV := V_\Z \otimes_{\Z_{\Cc}} \cO_\Cc.$$

Let us denote by $\mathrm{conj}$ the complex conjugation on the vector bundle $\cV$ induced by the above integral structure on $\cV$, and for every integer $r$, let us denote by $\ccH^{r,n-r}$ the real analytic vector subbundle of $\cV$ defined by:
$$\ccH^{r,n-r} := \cF^r \cap \mathrm{conj}(\cF^{n-r}).$$
Observe that the following equality of real analytic vector bundles holds:
$$\cV = \bigoplus_{r = 0}^n \ccH^{r,n-r}.$$

Let us assume that the variation of Hodge structures $\V$ is \emph{polarized}, namely that there exists a real bilinear form $Q$ on the real local system $V_\R := V_\Z  \otimes_\Z \R$, naturally extending into an $\cO_\Cc$-bilinear form on the vector bundle $\cV = V_\R \otimes_{\R_{\Cc}} \cO_\Cc$ (also denoted by $Q$), which satisfies the three following conditions:
\begin{itemize}
    \item the form $Q$ is symmetric if $n$ is even and antisymmetric if $n$ is odd,
    \item for every pair of integers $(r,r')$ such that $r \neq r'$, we have $Q(\ccH^{r,n-r}, \ccH^{n-r',r'}) = 0,$
    \item the following hermitian form on $\cV$:
    $$\bigg (u = \sum_{r = 0}^n u_r, v \bigg ) \longmapsto \sum_{r = 0}^n \Big (\sqrt{-1} \Big)^{2r-n} Q ( u_r, \mathrm{conj}(v)),$$
    where for every $r$, $u_r$ denotes an element of $\ccH^{r,n-r}$, is positive definite.
\end{itemize}

Such a polarization $Q$ induces, for every integer $p$, a hermitian metric $h_p$ on the vector bundle~$\cF^p/\cF^{p+1}$, and therefore a hermitian metric $h$ on the line bundle on $\Cc$:
$$\GK_\Cc(\V) = \bigotimes_{p \geq 0} (\det \cF^p / \cF^{p+1})^{\otimes p}.$$

The Chern form of the hermitian line bundle $\GK_\Cc(\V)$ is defined as the $(1,1)$-form on $\Cc$:
$$c(\GK_\Cc(\V), h) := (2\pi\sqrt{-1})^{-1} \, \partial \overline{\partial} \log h(s,s),$$
where $s$ denotes any local non-vanishing analytic section of $\GK_\Cc(\V)$. 
It is a $C^\infty$ real form on $\Cc$, and as shown by Griffiths \cite[Theorem 5.2 and Proposition 7.15]{Griffiths70}, it is non-negative.\footnote{This positivity result holds more generally for the real $(1,1)$-form  $c(\GK_S(\V), h)$ attached to some polarized variation of Hodge structures $\V$ over an arbitrary complex analytic manifold $S$.}

\begin{theorem}[{\cite[Proposition (3.4)]{Peters84}}] \label{GK curvature uni}
With the above notation, if the local monodromy of the variation of Hodge structures at every point of $\Delta$ is unipotent, then the non-negative $(1,1)$-form $c(\GK_\Cc(\V), h)$ is integrable on $\Cc$, and the following equality of integers holds:
\begin{equation} \label{eq GK curvature uni}
\hgt_{GK}(\V_\eta) = \int_{\Cc} c(\GK_\Cc(\V), h).
\end{equation}
\end{theorem}

Peters' proof crucially relies on a description of  the asymptotic behavior 
of the metric $h$ near the points of $\Delta = C \smallsetminus \Cc$, which follows from  Schmid's estimates \cite{Schmid73} concerning degenerations of Hodge structures. 

The following result is a straightforward consequence of  the definition of the stable Griffiths height, of the functoriality of the Chern form $c(\GK_\Cc(\V), h)$ under pull-back, and of Theorem \ref{GK curvature uni} applied to $C'\smallsetminus \sigma^{-1}(\Delta)$, where $\sigma: C' \ra C$ denotes a finite covering of $C$ such that the variation of Hodge structures $\sigma^\ast \V$ on $C'\smallsetminus \sigma^{-1}(\Delta)$ admits unipotent local monodromy at every point of $\sigma^{-1}(\Delta)$. 

\begin{corollary}\label{GK curvature stab} With the above notation, the non-negative $(1,1)$-form $c(\GK_\Cc(\V), h)$ is integrable on~$\Cc$,
and the following equality of rational numbers holds:
\begin{equation} \label{eq GK curvature stab}
\hgt_{GK, stab}(\V_\eta) = \int_{\Cc} c(\GK_\Cc(\V), h).
\end{equation}
\end{corollary}

\subsection{A semicontinuity theorem} Now we are in a position to state and to establish the semicontinuity property of the stable Griffiths height.

Let $D = D(0,1)$ be the open unit complex disk, let $S$ be an analytic surface, and let
$$p : S \lra D$$
be a proper smooth complex morphism. We shall assume that $S$, or equivalently every fiber of $p$, is connected.

Let $\Delta$ be a reduced analytic divisor in $S$ which is horizontal, namely which does not contain any fiber of $p$, or equivalently such that the restriction $p_{\mid \Delta}$ is a finite flat morphism. 

Let moreover $\V$ be an integral  polarized variation of Hodge structures on $S \smallsetminus \Delta$.

For every point $z$ in $D$, the fiber $S_z:= p^{-1}(z)$ is a connected smooth projective complex curve, and we shall denote its generic point by $\eta_z$.
Moreover the fiber $$\Delta_z:= p_{\mid \Delta}^{-1}(z) = \Delta \cap S_z$$ of $p_{\mid \Delta}$ is a finite subset in $S_z$, and we may consider the stable Griffiths height  $\hgt_{GK, stab}(\V_{\eta_z})$ of the  integral polarized variation of Hodge structures $\V_{\mid S_z \smallsetminus \Delta_z}$ on the Zariski open subset  $S_z \smallsetminus \Delta_z$ of the connected smooth projective complex curve $S_z$.

\begin{theorem}\label{semicontinuity GK}  
For every $z_0 \in D$, there exists an open neighborhood $\Omega$ of $z_0$ in $D$ and $\delta \in \Q$ such that:
  \begin{equation}\label{eq semicontinuity GK}
  \hgt_{GK, stab}(\V_{\eta_z}) = \delta \quad \mbox{for every $z \in  \Omega \smallsetminus \{z_0\}$},
   \end{equation} 
and:
    \begin{equation}\label{ineq semicontinuity GK}
\hgt_{GK, stab}(\V_{\eta_{z_0}}) \leq \delta.
    \end{equation}
\end{theorem}

\subsection{Proof of the semicontinuity theorem \ref{semicontinuity GK}}

 Since the reduced divisor $\Delta$ is horizontal relatively to $p$, the point $z_0$ admits a connected open neighborhood  $\Omega$ in $D$ such that, if we let
 $$\Omegac :=  \Omega \smallsetminus \{z_0\},$$ then the finite morphism $$p_{\mid \Delta \cap p^{-1}(\Omegac)}: \Delta \cap p^{-1}(\Omegac) \lra \Omegac$$
 is \'etale. Let us show that $\hgt_{GK, stab}(\V_{\eta_z})$ is a locally constant (and consequently constant) function of $z\in \Omegac$.
 
 To achieve this, observe that the conjugacy class of the (quasi-unipotent) local monodromy of the VHS $\V_{\mid S_{p(\omega)} \smallsetminus \Delta_{p(\omega)}}$ at some point $\omega \in \Delta \cap p^{-1}(\Omegac)$ is locally constant on $\Delta \cap p^{-1}(\Omegac)$.
 
 Assume  that there exists  some point $z_1$ of $\Omegac$ such that the local monodromy of $\V_{\mid S_{z_1} \smallsetminus \Delta_{z_1}}$ at every point $\omega \in \Delta_{z_1}$ is unipotent. Then the local monodromy of the VHS $\V_{\mid p^{-1}(\Omegac) \smallsetminus \Delta}$ along $\Delta \cap p^{-1}(\Omegac)$ is unipotent, and Peters' construction of the Griffiths line bundle makes sense\footnote{This is indeed the case for Deligne's construction \cite{Deligne70, Katz76}, and for Schmid's results on extensions of Hodge bundles \cite{Schmid73}.} for this VHS and defines some holomorphic line bundle $\GK_{p^{-1}(\Omegac)}(\V_{\mid p^{-1}(\Omegac) \smallsetminus \Delta})$ over $p^{-1}(\Omegac)$ whose restriction to $S_z$ is canonically isomorphic to $\GK_{S_z}(\V_{\eta_z})$ for every $z\in \Omegac$. Consequently the degree 
 $$\hgt_{GK}(\V_{\eta_z}) = \hgt_{GK, stab}(\V_{\eta_z})$$
 of $\GK_{S_z}(\V_{\eta_z})$ is independent of the choice of $z \in \Omegac$.

 In general, a point $z_1$ as above may not exist. However, over sufficiently small open subsets $\widetilde{\Omega}$ of~$\Omegac$, one reduces to this case by introducing a cyclic covering $S'_{\widetilde{\Omega}}$ of~$S_{\widetilde{\Omega}}:= p^{-1}(\widetilde{\Omega})$ suitably ramified over~$\Delta \cap S_{\widetilde{\Omega}}$ and smooth over~$\widetilde{\Omega}$.

  This establishes the existence of $\Omega$ and $\delta$ such that \eqref{eq semicontinuity GK} holds. 
  
  The inequality \eqref{ineq semicontinuity GK}, which indeed may be written as:
 $$\hgt_{GK, stab}(\V_{\eta_{z_0}}) \leq \liminf_{z \ra z_0} \hgt_{GK, stab}(\V_{\eta_z}),$$
 is an avatar of Fatou's lemma applied to the  positive $(1,1)$-form $c(\GK_{S\smallsetminus \Delta}(\V), h)$ attached to the polarized VHS $\V$ over $S\smallsetminus \Delta.$ Actually it will follow from the expression
 $$\hgt_{GK, stab}(\V_{\eta_z}) = \int_{S_z \setminus \Delta_z} c(\GK_{S\setminus \Delta}(\V), h)$$
 for the Griffiths heights $\hgt_{GK, stab}(\V_{\eta_z})$ established in Corollary \ref{GK curvature stab}, and from the positivity and continuity of the relative $(1,1)$-form defined by the restriction of  $c(\GK_{S\setminus \Delta}(\V), h)$ to the fibers of 
 $$p_{\mid S \setminus \Delta}: S\setminus \Delta \lra D.$$
 
 Indeed, after possibly shrinking the open neighborhood $\Omega$ of $z_0$ in $D$, we may find a neighborhood~$T$ of $\Delta_\Omega$ in $S_\Omega:= p^{-1}(\Omega)$, which is a closed $C^\infty$-submanifold with boundary of $S_\Omega$, whose boundary~$\partial T$ is transverse to the fibers of $p$. Moreover $T_{z_0} := T \cap S_{z_0}$ may be chosen to be an arbitrarily small neighborhood of $\Delta_{z_0}$ in $S_{z_0}$. 
 
 Then the complement in $S_\Omega$ of the interior $\mathaccent23{T}$ of $T$ defines  a $C^\infty$-manifold $S_\Omega \smallsetminus \mathaccent23{T}$ with boundary~$\partial T$ transverse to the fibers of $p$. The continuity of $c(\GK_{S\setminus \Delta}(\V), h)$ over $S\smallsetminus \Delta$, hence over $S_\Omega \smallsetminus \mathaccent23{T}$, shows that, if $\mathaccent23{T}_z := S_z \cap \mathaccent23{T},$ then the integral
 $$\int_{S_z \setminus \mathaccent23{T}_z} c(\GK_{S\setminus \Delta}(\V), h)$$
 is a continuous function of $z \in \Omega$. Moreover it satisfies, for every $z \in \Omega\smallsetminus \{z_0\},$ the inequality:
 $$\int_{S_z \setminus \mathaccent23{T}_z} c(\GK_{S\setminus \Delta}(\V), h) \leq \int_{S_z \setminus \Delta_z} c(\GK_{S\setminus \Delta}(\V), h) = \delta.$$
 
 Consequently:
$$\int_{S_{z_0} \setminus \mathaccent23{T}_{z_0}} c(\GK_{S\setminus \Delta}(\V), h) =
\lim_{z \ra z_0} \int_{S_z \setminus \mathaccent23{T}_z} c(\GK_{S\setminus \Delta}(\V), h)  \leq \delta.$$
Using Corollary  \ref{GK curvature stab}, and the fact that the neighborhood $T_{z_0}$ of $\Delta_{z_0}$ in $S_{z_0}$ can be chosen to be arbitrarily small, this implies the inequality: 
\begin{equation*}
\hgt_{GK, stab}(\V_{\eta_{z_0}}) = \int_{S_{z_0} \setminus \Delta_{z_0}} c(\GK_{S\setminus \Delta}(\V), h)   \leq \delta.  
\end{equation*}

\section{The inequality $\hgt_{GK, stab}\big(\H^{N-1}(H_\eta/C_\eta)\big)\leq F_{stab}(d,N) \,\hgt_{int}(H/C)$}
\label{sec: ineq GK semi-hom}

In this final section, we complete the proofs of our results concerning the validity of the inequality \eqref{ineq GK GIT semihom}  under suitable assumptions, 
alluded to in Paragraph \ref{1.7}. 

They  will follow from the estimate
$$\hgt_{GK, stab}\big(\H^{N-1}(H_\eta/C_\eta)\big)\leq F_{stab}(d,N) \,\hgt_{int}(H/C),$$
which we will establish when the hypersurface $H_\eta$ admits  a ``model'' $H$ over $C$ which is a pencil of hypersurfaces  in a projective bundle $\PP(E)$, and when moreover one of the following conditions is satisfied: $H$ is smooth and its fibers over $C$ admit  at most one singular point, assumed to be semihomogeneous (Proposition \ref{intro upper bound GK}); or  the instability $\mu_{\max}(E) -\mu(E)$ of $E$ is not too large when compared to $\hgt_{int}(H/C)$ (Proposition \ref{generization}).

\subsection{The quotient $w_{N,\delta}(\delta-1)^{-N}$ is a non-increasing function of $\delta$}

Recall that, in Theorem ~\ref{intro GK hyp P(E) hom crit},  to any pair $(N, \delta)$ of positive integers, we have attached the rational number 
$$w_{N,\delta} = (\delta-1)\big[(N \delta + 1) (\delta - 1)^{N-1} + (-1)^N (\delta + 1) \big]/(12 \delta^2).$$

\begin{proposition}\label{decreasing contribution}
Let $N$ be a positive integer.
The rational-valued sequence $(g_N(\delta))_{\delta \in \N_{\geq 2}}$ defined by:
\begin{equation}\label{defg}
g_N(\delta) :=12 \, w_{N,\delta}(\delta-1)^{-N} = \big[(N \delta + 1) (\delta-1)^{N-1} + (-1)^N (\delta + 1)\big] \, \delta^{-2} (\delta -1)^{-N +1}
\end{equation}
is non-increasing.
\end{proposition}

Observe that the sequence $g_1$ is the zero sequence. It will follow from the proof that the sequence $(g_N(\delta))_{\delta \in \N_{\geq 2}}$ is actually decreasing when $N=2$ or $N \geq 4$. Moreover $g_3(2) = g_3(3)= 1$ and the sequence $(g_3(\delta))_{\delta \in \N_{\geq 3}}$ is decreasing.

\begin{proof} The result is clear if $N = 1,$ so we shall assume $N\geq 2. $

The sequence $(g_N(\delta))_{\delta \in \N_{\geq 2}}$ can clearly be extended into a real-valued function on $[2, +\infty)$  still given  by formula \eqref{defg}. We shall also denote this function by $g_N.$

A straightforward computation yields that the derivative of $g_N$ is given by:
\begin{equation} \label{eq deriv g}
\frac{d g_N(\delta)}{d\delta} = -1/\big(\delta^3 (\delta-1)^N\big) \, [(N \delta + 2)  (\delta-1)^N + (-1)^N (N \delta^2 + N \delta - 2) ].
\end{equation}

The terms $N \delta + 2$ and $N \delta^2 + N \delta - 2$ are positive for $N \geq 2$ and $\delta \geq 2$, so the numerator in \eqref{eq deriv g} admits the following lower bound:
\begin{align*}
(N \delta + 2)  (\delta-1)^N + (-1)^N (N \delta^2 + N \delta - 2)
& \geq (N \delta + 2) (\delta-1)^N - (N \delta^2 + N \delta - 2) \\
& \geq (N \delta + 2) (\delta-1)^2 - (N \delta^2 + N \delta - 2) \\
& = N \delta^2 (\delta-3) + 2 \delta (\delta-2) + 4,
\end{align*}
which is positive when $\delta \geq 3.$ Consequently the derivative $d g_N(\delta)/d\delta$ is negative for $\delta \geq 3$, and the function $g_N$ is  decreasing on $[3, +\infty).$  In particular, the sequence $(g_N(\delta))_{\delta \in \N_{\geq 3}}$ is decreasing.

It remains to prove the inequality $g_N(2) \geq g_N(3).$ The values of $g_N$ at 2 and 3 are given by:
$$g_N(2) = (1 / 4) [2 N  + 1 + 3 (-1)^N ]$$
and:
$$g_N(3) = 1/(9 . 2^{N-1})  [(3 N + 1) 2^{N-1} + 4 (-1)^N],$$
and their difference is given by:
\begin{align}
g_N(2) - g_N(3) & = 1/(9 . 2^N)  [9 . 2^{N-2}  (2N + 1 + 3 (-1)^N) - (3 N + 1) 2^N - 8 (-1)^N ] \nonumber \\
& = 1/(9 . 2^N)   [2^{N-2} (18 N + 9 + 27 (-1)^N - 12 N - 4 ) - 8 (-1)^N] \nonumber \\
\label{eq diff 2 3} & = 1/(9 . 2^N)  [2^{N-2} (6 N + 5 + 27 (-1)^N) - 8 (-1)^N].
\end{align}

If $N$ is even, the numerator in \eqref{eq diff 2 3} is  $2^{N-2} (6 N + 32) - 8$ and is positive.
If $N$ is odd, the numerator in \eqref{eq diff 2 3} is  $2^{N-2} (6 N - 22) + 8$; it vanishes if $N = 3$ and is positive if $N\geq  5.$
Consequently the difference $g_N(2) - g_N(3)$ is always non-negative, as wanted.
\end{proof}

\subsection{Upper bound and equality cases for the stable Griffiths height of pencils of hypersurfaces with semihomogeneous singularities}

We are now in position to complete the proof of Proposition \ref{intro upper bound GK}. 

Recall that $F_{stab}(d,N)$ is the element in $(1/12) \Z$ given when $N$ is odd by:
\begin{align*}
F_{stab}(d,N) & := \frac{N+1}{24 d^2} \left[ (d-1)^N (d^2 N - d^2 - 2 d N - 2 ) + 2 (d^2-1) \right] \\
& = (N+1)(d-1)/(24 d^2) [(d-1)^{N-1} (d^2 N - d^2 - 2 d N - 2 ) + 2 (d+1)],
\end{align*}
and, when $N$ is even, by:
\begin{align*}
F_{stab}(d,N) &:= \frac{N+1}{24 d^2} \left [ (d-1)^N  (d^2 N + 2 d^2 - 2 d N - 2) - 2 (d^2-1) \right ] \\
& = (N+1)(d-1)/(24 d^2) \left[(d-1)^{N-1} (d^2 N + 2 d^2 - 2 d N - 2) - 2 (d+1) \right ].
\end{align*}

For every point $P$ in $\Sigma$,  the integer $\delta_P$ is at least 2, and therefore according to  Proposition \ref{decreasing contribution} the following inequality of rational numbers holds:
$$ 12\,  w_{N,\delta_P} (\delta_P-1)^{-N} = g_N(\delta_P) \leq g_N(2) = (1 / 4) [2 N  + 1 + 3 (-1)^N ].$$
We can rewrite it as follows:
\begin{equation}\label{upper bound v}
 w_{N,\delta_P} \leq  \big((\delta_P-1)^N / 48\big)   [2 N + 1 + 3 (-1)^N].
 \end{equation}

Consequently, applying equality \eqref{intro GKstabXL P(E) hom crit}, the following inequality holds:
\begin{align}
\hgt_{GK, stab}\big(\H^{N-1}(H_\eta/C_\eta)\big) 
& =  - (N+1) w_{N,d}  \; \hgt_{int}(H/C) +\sum_{P \in \Sigma} w_{N,\delta_P} \nonumber  \\
& \leq  - (N+1) w_{N,d} \; \hgt_{int}(H/C) + \sum_{P \in \Sigma} \big((\delta_P-1)^N / 48\big) [2 N + 1 + 3 (-1)^N]  \nonumber \\
\label{first eq upper bound stab}
&= - (N+1) w_{N,d} \; \hgt_{int}(H/C)  \\
  & \quad  + (1/48) \, [2 N + 1 + 3 (-1)^N]  (N+1) (d-1)^N \hgt_{int}(H/C) \nonumber \\
&=  F(d,N) \; \hgt_{int}(H/C),  \nonumber
\end{align}
where in \eqref{first eq upper bound stab}, we have applied equality \eqref{intro f(sigma) P(E) hom crit}, and where $F(d,N)$ is the rational number defined by:
\begin{align*}
F(d,N) &:= - (N+1) w_{N,d} + (N+1) (d-1)^N/48  \, [2 N + 1 + 3 (-1)^N]\\
&
= - (N+1) (d-1)/(12 d^2)\,  [(N d + 1) (d - 1)^{N-1} + (-1)^N (d + 1) ] \\
& \quad + (N+1) (d-1)^N/48 \, [2 N + 1 + 3 (-1)^N] \\
& = (N+1) (d-1)/(48 d^2) \, [- 4 (N d + 1) (d-1)^{N-1} - 4 (-1)^N (d+1) \\
& \quad + (d-1)^{N-1} d^2 (2 N + 1 + 3 (-1)^N) ] \\
& = (N+1) (d-1)/(48 d^2) \, [(d-1)^{N-1} (2 d^2 N + d^2 (1 + 3 (-1)^N) - 4 d N - 4) - 4 (-1)^N (d+1)].
\end{align*}

If $N$ is odd, this is given by:
\begin{align*}
F(d,N) & = (N+1) (d-1)/(48 d^2) \, [(d-1)^{N-1} (2 d^2 N - 2 d^2 - 4 d N - 4) + 4 (d+1)] \\
& = (N+1) (d-1)/(24 d^2) \, [(d-1)^{N-1} (d^2 N - d^2 - 2 d N - 2) + 2 (d+1)] \\
& = F_{stab}(d,N),
\end{align*}
and if $N$ is even, by:
\begin{align*}
F(d,N) &= (N+1) (d-1)/(48 d^2) \, [(d-1)^{N-1} (2 d^2 N + 4 d^2 - 4 d N - 4) - 4 (d+1)] \\
& = (N+1) (d-1)/(24 d^2) \, [(d-1)^{N-1} (d^2 N + 2 d^2 - 2 d N - 2) - 2 (d+1)] \\
& = F_{stab}(d,N).
\end{align*}

In both cases, $F(d,N)$ is precisely $F_{stab}(d,N)$, as wanted.

The assertions in Proposition \ref{intro upper bound GK} concerning the equality cases follow from the above proof and from the observations concerning the decreasing character of $(g_N(\delta))_{\delta \in \N_{\geq 2}}$ stated after Proposition~\ref{decreasing contribution}.

\subsection{Lower semicontinuity of $\hgt_{GK, stab}$ and specialization of pencils of hypersurfaces with generic singularities}\label{proof generization GIT}

\subsubsection{} In this subsection, we will establish the estimate \eqref{pref intro ineq GK int hom Bis} in Proposition \ref{generization} as a consequence of Theorem \ref{pref intro GK hyp P(E)}, of the semicontinuity theorem (Theorem \ref{semicontinuity GK}), and of the following genericity result. 

\begin{theorem}[{\cite[Th\'eor\`eme A.1.3]{MordantMem23}}]\label{th genericity in GIT}
Let $k$ be an algebraically closed field of characteristic zero, $C$ a connected smooth projective $k$-curve of genus $g$, $E$ a non-zero vector bundle over $C$, $M$ a line bundle over $C$, and $d \geq 1$ an integer. Let moreover:
$$\pi: \PP(E) := \mathrm{Proj}_C(S^\bullet E^\vee) \lra C$$
be the projective bundle associated to $E$.

If the degree of $M$ satisfies the following lower bound:
$$\deg_C M > 2g-1+d \mu_{\max}(E),$$
then there exists a non-empty Zariski open subset $U$ of the vector space of sections $H^0(\PP(E),\cO_E(d) \otimes \pi^\ast M)$
such that, for every section $s$ in $U$, the hypersurface $H_s$ of $\PP(E)$ defined by the vanishing of~$s$ is smooth over $k$, every critical point of the restriction
$$\pi_{\mid H_s} : H_s \lra C$$
is non-degenerate, and every fiber of $\pi_{\mid H_s}$ contains at most one critical point.
\end{theorem}

We now adopt the notation and hypotheses of the first part of Proposition \ref{generization}.

Actually we may assume that $\hgt_{int}(H/C)$ satisfies the following stronger lower bound:
\begin{equation} \label{ht int bound 2 g - 1}
\hgt_{int}(H/C)  > 2g -1 + d(\mu_{\max}(E) -\mu(E)).
\end{equation}
Indeed we are reduced to the case where this inequality holds by replacing $C$, supposed to have positive genus $g$,  with a connected unramified covering $C'$ of $C$ of sufficiently large degree $\delta$, using the compatibility with finite base change of the slopes $\mu(E)$ and $\mu_{\max}(E)$ and of the heights $\hgt_{int}(H/C)$ and  $\hgt_{GK, stab}\big(\H^{N-1}(H_\eta/C_\eta)\big)$.\footnote{If $g'$ denotes the genus of $C'$, according to the Riemann-Hurwitz formula, we have: $(2g'-1)/\delta = 2g -2 + 1/\delta$, and this ratio converges to $2g-2$ when $\delta$ goes to infinity.}

\subsubsection{} The structure of the Chow group of the projective bundle $\PP(E)$ (see for instance \cite[Th. 3.3 (b)]{Fulton98}) implies that the divisor $H$ in $\PP(E)$ is defined by the vanishing of a section
$$s \in H^0(\PP(E),\cO_E(d) \otimes \pi^\ast M)$$
where $M$ is some line bundle over $C$.

Recall, with this notation, the alternate expression of $\hgt_{int}(H/C)$ given in \cite[(6.2.13)]{Mordant22}:
\begin{equation} \label{htint M E in GIT deformation}
\hgt_{int}(H/C) = \deg_C M - \frac{d}{N+1} \deg_C E.
\end{equation}
Using 
equality \eqref{htint M E in GIT deformation} and inequality \eqref{ht int bound 2 g - 1}, we obtain the following lower bound for $\deg_C M$:
$$\deg_C M = \hgt_{int}(H/C) + d  \mu(E) > 2g - 1 + d \mu_{\max}(E).$$

Consequently we may consider a non-empty Zariski open subset $U$ in $H^0(\PP(E),\cO_E(d) \otimes \pi^\ast M)$ as in Theorem \ref{th genericity in GIT} (applied here for $k := \C$). Since the vector space $H^0(\PP(E),\cO_E(d) \otimes \pi^\ast M)$ is irreducible for the Zariski topology, the subset $U$ is dense for the Zariski topology, hence also for the analytic topology. If $D = D(0,1)$ denotes the open unit complex disk, we may therefore consider a complex analytic map:
$$\phi : D \lra H^0(\PP(E),\cO_E(d) \otimes \pi^\ast M)$$
which maps $0$ to the section $s$ and maps the complement:
$$D^\ast := D \setminus \{0\}$$
into the open subset $U$. 

Let us define a complex analytic surface by:
$$S := C \times D,$$
let us denote the two projections by:
$$\mathrm{pr}_1 : C \times D \lra C \quad \mbox{and} \quad \mathrm{pr}_2 : C \times D \lra D,$$
and let us define a divisor in $\PP(E) \times D$ by:
$$H_\phi := \{(x,z) \in \PP(E) \times  D \mid \phi(z)(x) = 0\}.$$

For every $z$ in $D$, let us denote by $H_{\phi, z}$  the inverse image of the closed analytic subspace $C \times \{z\}$ in $S$ by the morphism:
$$(\pi \times \mathrm{Id}_D)_{\mid H_\phi} : H_\phi \subset \PP(E) \times D \lra S := C \times D.$$
It is precisely the divisor in the projective bundle $\PP(E) \times \{z\} \simeq \PP(E)$ over the curve $C \times \{z\} \simeq C$ defined by the vanishing of the section $\phi(z)$ in $H^0(\PP(E),\cO_E(d) \otimes \pi^\ast M)$. In particular, for $z = 0$, since $\phi(0) = s$, the divisor $H_{\phi, 0}$ may be identified with the pencil of hypersurfaces $H$ over $C$.

Observe moreover that for every $z$ in $D$, the divisor $H_{\phi, z}$ is horizontal over $C$ (hence a pencil of projective hypersurfaces) and that its generic fiber, denoted by $H_{\phi, (\eta,z)}$, is non-singular. Indeed, if~$z = 0$ this follows from the above identification and the hypothesis on the generic fiber $H_\eta$, and if~$z \neq 0$ this follows from the fact that $\phi(z)$ is in $U$ and from the construction of $U$.

Consequently there exists a reduced analytic divisor $\Delta$ in $S$, which is horizontal relatively to the proper smooth complex morphism
$$p := \mathrm{pr}_2 : S := C \times D \lra D,$$
and satisfies that the following proper complex morphism:
$$(\pi \times \mathrm{Id}_D)_{\mid H_\phi} : H_\phi \subset \PP(E) \times D \lra S := C \times D$$
is smooth over $S \smallsetminus \Delta$. 

Let us consider the inverse image:
$$H_{\phi, S \setminus \Delta} :=(\pi \times \mathrm{Id}_D)_{\mid H_\phi}^{-1}(S \setminus \Delta), $$
and the variation of Hodge structures $\V$ on $S \smallsetminus \Delta$ defined as the relative middle-dimensional cohomology of $H_{\phi, S \smallsetminus \Delta}$ over $S \smallsetminus \Delta$:
$$\V := \H^{N-1}(H_{\phi, S \setminus \Delta} / (S \setminus \Delta)).$$

\subsubsection{} Applying Theorem \ref{semicontinuity GK} to this variation of Hodge structures, we obtain a neighborhood $\Omega$ of~$0$ in $D$ and a rational number~$\delta$ such that:
\begin{equation}\label{eq semicontinuity GK for GIT}
  \hgt_{GK, stab}\big(\H^{N-1}(H_{\phi,(\eta, z)} / (C \times \{z\})_{(\eta, z)})\big) = \delta \quad \mbox{for every $z \in  \Omega \smallsetminus \{0\}$},
   \end{equation} 
and:
    \begin{equation}\label{ineq semicontinuity GK for GIT}
\hgt_{GK, stab}\big(\H^{N-1}(H_{\phi,(\eta, 0)}/ (C \times \{0\})_{(\eta, 0)})\big) \leq \delta.
    \end{equation}

Let us fix an element $z$ in $\Omega \smallsetminus \{0\}$. By construction, the pencil of hypersurfaces $H_{\phi,z}$ over $C \times \{z\} \simeq C$ satisfies the smoothness and non-degeneracy hypotheses of Theorem \ref{pref intro GK hyp P(E)}, so that the following equality of rational numbers holds:
\begin{equation} \label{eq semicontinuity generic sing} 
\hgt_{GK, stab}\big(\H^{N-1}(H_{\phi,(\eta, z)}/ (C \times \{z\})_{(\eta, z)})\big) = F_{stab}(d,N) \, \hgt_{int}(H_{\phi, z}/C).
\end{equation}

As the divisor $H_{\phi, z}$ (resp. $H$) in $\PP(E) \times \{z\} \simeq \PP(E)$ (resp. $\PP(E)$) is defined by the vanishing of the section $\phi(z)$ (resp. $s$) in $H^0(\PP(E),\cO_E(d) \otimes \pi^\ast M)$, the following equality of rational numbers holds as a consequence of equality \eqref{htint M E in GIT deformation} and the analogue of this equality for the pencil $H_{\phi, z}$:
\begin{equation}\label{eq contin ht int}
\hgt_{int}(H_{\phi, z}/C) 
=\hgt_{int}(H/C).
\end{equation}

The desired inequality \eqref{pref intro ineq GK int hom Bis} is a straightforward consequence of inequality \eqref{ineq semicontinuity GK for GIT}, of the identification between the pencil of hypersurfaces $H_{\phi, 0}$ over $C \times \{0\} \simeq C$ with the pencil $H$ over $C$, and of equalities \eqref{eq semicontinuity GK for GIT} (applied to the chosen $z$), \eqref{eq semicontinuity generic sing}, and \eqref{eq contin ht int}.


\begin{thebibliography}{AHLH23}

\bibitem[AHLH23]{AHLH23}
J.~Alper, D.~Halpern-Leistner, and J.~Heinloth.
\newblock Existence of moduli spaces for algebraic stacks.
\newblock {\em Invent. Math.}, 234(3):949--1038, 2023.

\bibitem[Alp13]{Alper13}
J.~Alper.
\newblock Good moduli spaces for {A}rtin stacks.
\newblock {\em Ann. Inst. Fourier (Grenoble)}, 63(6):2349--2402, 2013.

\bibitem[Bea09]{Beauville09}
A.~Beauville.
\newblock Moduli of cubic surfaces and {H}odge theory (after {A}llcock,
  {C}arlson, {T}oledo).
\newblock In {\em G\'{e}om\'{e}tries \`a courbure n\'{e}gative ou nulle,
  groupes discrets et rigidit\'{e}s}, volume~18 of {\em S\'{e}min. Congr.},
  pages 445--466. Soc. Math. France, Paris, 2009.

\bibitem[Ben14]{Benoist14}
O.~Benoist.
\newblock Quelques espaces de modules d'intersections compl\`etes lisses qui
  sont quasi-projectifs.
\newblock {\em J. Eur. Math. Soc. (JEMS)}, 16(8):1749--1774, 2014.

\bibitem[Bog78]{Bogomolov78}
F.~A. Bogomolov.
\newblock Holomorphic tensors and vector bundles on projective manifolds.
\newblock {\em Izv. Akad. Nauk SSSR Ser. Mat.}, 42(6):1227--1287, 1439, 1978.

\bibitem[Bos94]{Bost94}
J.-B. Bost.
\newblock Semi-stability and heights of cycles.
\newblock {\em Invent. Math.}, 118(2):223--253, 1994.

\bibitem[Bos96]{Bost96}
J.-B. Bost.
\newblock Intrinsic heights of stable varieties and abelian varieties.
\newblock {\em Duke Math. J.}, 82(1):21--70, 1996.

\bibitem[Del70]{Deligne70}
P.~Deligne.
\newblock {\em \'{E}quations diff\'{e}rentielles \`a points singuliers
  r\'{e}guliers}.
\newblock Lecture Notes in Mathematics, Vol. 163. Springer-Verlag, Berlin-New
  York, 1970.

\bibitem[Dem12]{Demazure12}
M.~Demazure.
\newblock R\'{e}sultant, discriminant.
\newblock {\em Enseign. Math. (2)}, 58(3-4):333--373, 2012.

\bibitem[Fal83]{Faltings83}
G.~Faltings.
\newblock Endlichkeitss\"{a}tze f\"{u}r abelsche {V}ariet\"{a}ten \"{u}ber
  {Z}ahlk\"{o}rpern.
\newblock {\em Invent. Math.}, 73(3):349--366, 1983.

\bibitem[FC90]{Faltings-Chai90}
G.~Faltings and C.-L. Chai.
\newblock {\em Degeneration of abelian varieties}, volume~22 of {\em Ergebnisse
  der Mathematik und ihrer Grenzgebiete (3)}.
\newblock Springer-Verlag, Berlin, 1990.
\newblock With an appendix by David Mumford.

\bibitem[Ful98]{Fulton98}
W.~Fulton.
\newblock {\em Intersection theory}, volume~2 of {\em Ergebnisse der Mathematik
  und ihrer Grenzgebiete (3)}.
\newblock Springer-Verlag, Berlin, second edition, 1998.

\bibitem[Gas00]{Gasbarri00}
C.~Gasbarri.
\newblock Heights and geometric invariant theory.
\newblock {\em Forum Math.}, 12(2):135--153, 2000.

\bibitem[GKZ08]{GKZ08}
I.~M. Gelfand, M.~M. Kapranov, and A.~V. Zelevinsky.
\newblock {\em Discriminants, resultants and multidimensional determinants}.
\newblock Modern Birkh\"{a}user Classics. Birkh\"{a}user Boston, Inc., Boston,
  MA, 2008.
\newblock Reprint of the 1994 edition.

\bibitem[Gri70]{Griffiths70}
P.~A. Griffiths.
\newblock Periods of integrals on algebraic manifolds. {III}. {S}ome global
  differential-geometric properties of the period mapping.
\newblock {\em Inst. Hautes \'{E}tudes Sci. Publ. Math.}, 38:125--180, 1970.

\bibitem[Hil90]{Hilbert90}
D.~Hilbert.
\newblock {\"U}ber die {Theorie} der algebraischen {Formen}.
\newblock {\em Math. Ann.}, 36:473--534, 1890.

\bibitem[Hil93]{Hilbert93}
D.~Hilbert.
\newblock Ueber die vollen {I}nvariantensysteme.
\newblock {\em Math. Ann.}, 42(3):313--373, 1893.

\bibitem[Jor80]{Jordan80}
C.~Jordan.
\newblock M\'{e}moire sur l'\'{e}quivalence des formes.
\newblock {\em J. \'{Ecole} Polytechnique}, 48:112--150, 1880.

\bibitem[Kat76]{Katz76}
N.~M. Katz.
\newblock An overview of {D}eligne's work on {H}ilbert's twenty-first problem.
\newblock In {\em Mathematical developments arising from {H}ilbert problems
  ({P}roc. {S}ympos. {P}ure {M}ath., {V}ol. {XXVIII}, {N}orthern {I}llinois
  {U}niv., {D}e {K}alb, {I}ll., 1974)}, pages 537--557, 1976.

\bibitem[Kat13]{Kato13}
K.~Kato.
\newblock Heights of mixed motives.
\newblock \texttt{https://arxiv.org/abs/1306.5693}, 2013.

\bibitem[Kat14]{Kato14}
K.~Kato.
\newblock Heights of motives.
\newblock {\em Proc. Japan Acad. Ser. A Math. Sci.}, 90(3):49--53, 2014.

\bibitem[Kat18]{Kato18}
K.~Kato.
\newblock Height functions for motives.
\newblock {\em Selecta Math. (N.S.)}, 24(1):403--472, 2018.

\bibitem[Kat20]{Kato20}
K.~Kato.
\newblock Height functions for motives, {II}.
\newblock In {\em Development of {I}wasawa theory---the centennial of {K}.
  {I}wasawa's birth}, volume~86 of {\em Adv. Stud. Pure Math.}, pages 467--535.
  Math. Soc. Japan, Tokyo, [2020] \copyright 2020.

\bibitem[Kos15]{Koshikawa15}
T.~Koshikawa.
\newblock On heights of motives with semistable reduction.
\newblock \texttt{https://arxiv.org/abs/1505.01873}, 2015.

\bibitem[Lee08]{Lee08}
Y.~Lee.
\newblock Chow stability criterion in terms of log canonical threshold.
\newblock {\em J. Korean Math. Soc.}, 45(2):467--477, 2008.

\bibitem[Mac17]{Maculan17}
M.~Maculan.
\newblock Diophantine applications of geometric invariant theory.
\newblock {\em M\'{e}m. Soc. Math. Fr. (N.S.)}, 152:1--149, 2017.


\bibitem[Mor22]{Mordant22}
T.~Mordant.
\newblock Griffiths heights and pencils of hypersurfaces.
\newblock \texttt{https://arxiv.org/abs/2212.11019}, to appear in
  \emph{M{\'e}m. Soc. Math. Fr.}, 2022.

\bibitem[Mor23]{MordantMem23}
T.~Mordant.
\newblock Hauteurs de {G}riffiths-{K}ato des pinceaux de vari\'et\'es
  projectives.
\newblock PhD dissertation, available at
  \texttt{https://www.imo.universite-paris-saclay.fr/en/perso/thomas-mordant/},
  2023.

\bibitem[Mor24]{Mordant23}
T.~Mordant.
\newblock A note on the semistability of singular projective hypersurfaces.
\newblock {\em Math. Z.}, 306(4):Paper No. 67, 19, 2024.

\bibitem[Mora]{MordantGIT2}
T.~Mordant.
\newblock Pencils of projective hypersurfaces, {G}riffiths heights and
  geometric invariant theory~{II}: Hypersurfaces with semihomogeneous
  singularities.
\newblock In preparation.


\bibitem[Morb]{MordantGMS}
T.~Mordant.
\newblock Good moduli spaces, heights, and finite semistable replacement.
\newblock {I}n preparation.



\bibitem[MB85a]{Moret-Bailly85a}
L.~Moret-Bailly.
\newblock Compactifications, hauteurs et finitude.
\newblock {\em Ast\'{e}risque}, 127:113--129, 1985.

\bibitem[MB85b]{Moret-Bailly85b}
L.~Moret-Bailly.
\newblock Pinceaux de vari\'{e}t\'{e}s ab\'{e}liennes.
\newblock {\em Ast\'{e}risque}, 129:1--266, 1985.

\bibitem[MFK94]{MumfordFogartyKirwan94}
D.~Mumford, J.~Fogarty, and F.~Kirwan.
\newblock {\em Geometric invariant theory}, volume~34 of {\em Ergebnisse der
  Mathematik und ihrer Grenzgebiete (2)}.
\newblock Springer-Verlag, Berlin, third edition, 1994.

\bibitem[Pet84]{Peters84}
C.~A.~M. Peters.
\newblock A criterion for flatness of {H}odge bundles over curves and geometric
  applications.
\newblock {\em Math. Ann.}, 268(1):1--19, 1984.

\bibitem[Sch73]{Schmid73}
W.~Schmid.
\newblock Variation of {H}odge structure: the singularities of the period
  mapping.
\newblock {\em Invent. Math.}, 22:211--319, 1973.

\bibitem[Ste77]{Steenbrink77}
J.~H.~M. Steenbrink.
\newblock Mixed {H}odge structure on the vanishing cohomology.
\newblock In {\em Real and complex singularities ({P}roc. {N}inth {N}ordic
  {S}ummer {S}chool/{NAVF} {S}ympos. {M}ath., {O}slo, 1976)}, pages 525--563,
  1977.

\bibitem[Ste76]{Steenbrink76}
J.~H.~M. Steenbrink.
\newblock Limits of {H}odge structures.
\newblock {\em Invent. Math.}, 31(3):229--257, 1975/76.

\bibitem[Szp85]{Szpiro85}
L.~Szpiro, editor.
\newblock {\em S\'{e}minaire sur les pinceaux arithm\'{e}tiques: la conjecture
  de {M}ordell}.
\newblock Soci\'{e}t\'{e} Math\'{e}matique de France, Paris, 1985.
\newblock Papers from the seminar held at the \'{E}cole Normale Sup\'{e}rieure,
  Paris, 1983--84, Ast\'{e}risque No. 127 (1985).

\bibitem[Tat66]{Tate66}
J.~Tate.
\newblock Endomorphisms of abelian varieties over finite fields.
\newblock {\em Invent. Math.}, 2:134--144, 1966.

\bibitem[Voi03]{Voisin03}
C.~Voisin.
\newblock {\em Hodge theory and complex algebraic geometry. {II}}, volume~77 of
  {\em Cambridge Studies in Advanced Mathematics}.
\newblock Cambridge University Press, Cambridge, 2003.

\bibitem[Voi22]{Voisin22}
C.~Voisin.
\newblock Schiffer variations and the generic {T}orelli theorem for
  hypersurfaces.
\newblock {\em Compos. Math.}, 158(1):89--122, 2022.

\bibitem[Zar75]{Zarhin75}
Ju.~G. Zarhin.
\newblock Endomorphisms of {A}belian varieties over fields of finite
  characteristic.
\newblock {\em Izv. Akad. Nauk SSSR Ser. Mat.}, 39(2):272--277, 471, 1975.

\bibitem[Zha96]{Zhang96}
S.~Zhang.
\newblock Heights and reductions of semi-stable varieties.
\newblock {\em Compositio Math.}, 104(1):77--105, 1996.

\end{thebibliography}

\end{document}